\pdfoutput=1       

\documentclass{llncs}
\usepackage{amsmath,amsfonts,amssymb}
\usepackage{graphicx,rotating}
\usepackage{color}

\newcommand{\llong}[1]{}
\newcommand{\short}[1]{#1}
\renewcommand{\llong}[1]{#1}
\renewcommand{\short}[1]{}

\newcommand{\anonymize}[1]{#1}
\newcommand{\anonymRepl}[2]{#1}  

\usepackage[citecolor=blue, pdftitle={\llong{Exploration of models}\short{Local search} for a cargo assembly planning problem}, pdfauthor={Belov,Boland,Savelsbergh,Stuckey}, pdfsubject={discrete optimization,constraint programming,cutting,packing,scheduling}, pdfkeywords={packing,cutting,scheduling,resource-constrained scheduling,resource constraint,non-preemptiveness,iterative heuristic,large neighbourhood search}]{hyperref}
\usepackage[comma,square,sort,numbers]{natbib}

\DeclareMathOperator{\veta}{eta}
\DeclareMathOperator{\tonn}{tonn}
\DeclareMathOperator{\speed}{speed}
\DeclareMathOperator{\rnd}{rnd}

\DeclareMathOperator{\cumulative}{\texttt{cumulative}}
\DeclareMathOperator{\disjunctive}{\texttt{disjunctive}}

\DeclareMathOperator{\delay}{delay}
\DeclareMathOperator{\depEarliest}{depEarliest}
\DeclareMathOperator{\card}{card}
\DeclareMathOperator{\stkDays}{stkDays}

\newcommand{\didpjs}{}
\newcommand{\NEW}{}

\newcommand{\GLEB}[1]{
#1}

\llong{\addtolength{\textheight}{4cm}
\addtolength{\topmargin}{-2cm}
\addtolength{\oddsidemargin}{-1cm}
\addtolength{\evensidemargin}{-1cm}
\addtolength{\textwidth}{2cm}}

\sloppypar

\newcommand{\Z}{{\mathbb{Z}}}

\newcommand{\be}{\begin{equation}}
\newcommand{\ee}{\end{equation}}

\begin{document}
\title{\llong{Exploration of models}\short{Local search} for a \\ cargo assembly planning problem}
\anonymRepl{\author{G.~Belov\inst{1} and N.~Boland\inst{2} and
  M.W.P.~Savelsbergh\inst{2} and P.J.~Stuckey\inst{1}}}{\author{}}
\anonymRepl{
\institute{Department of Computing and Information Systems \\
University of Melbourne, 3010 Australia
\llong{\\
E-mail: pstuckey@unimelb.edu.au}
\and
School of Mathematical and Physical Sciences \\
University of Newcastle, Callaghan 2308, Australia.
}}{\institute{}}

\date{\footnotesize Technical Report, 20 November 2014}
\maketitle
\begin{center}
{\footnotesize Technical Report, 20 November 2014}
\end{center}

\begin{abstract}
We consider a real-world cargo assembly planning problem arising in a coal
supply chain. The cargoes are built on the stockyard at a port terminal from
coal delivered by trains. Then the cargoes are loaded onto vessels. Only a
limited number of arriving vessels is known in advance. The goal is to
minimize the average delay time of the vessels over a long planning
period. We model the problem in the MiniZinc constraint programming language
and design a large neighbourhood search scheme.
\llong{
The effects of various optional constraints are investigated. \GLEB{Some of the optional constraints expand the model's scope toward a system view (berth capacity, port channel)}. An adaptive
scheme for a greedy heuristic from the literature is proposed and compared
to the constraint programming approach.}
\short{We compare against (an extended version of) a greedy heuristic
for the same problem.}
\end{abstract}

{\bf Keywords:} packing, scheduling, resource constraint, large neighbourhood search, constraint programming, adaptive greedy, visibility horizon

%
\section{Introduction}
%


The Hunter Valley Coal Chain (HVCC) refers to the inland portion of
the coal export supply chain in the Hunter Valley, New South Wales,
Australia. \llong{ Most of the coal mines in the Hunter Valley are
open pit mines. The coal is mined and stored either at a railway
siding located at the mine or at a coal loading facility used by
several mines. The coal is then transported to one of the terminals
at the Port of Newcastle, almost exclusively by rail. The coal is
dumped and stacked at a terminal to form stockpiles. } Coal from
different mines with different characteristics is `mixed' in a
stockpile \short{at a terminal at the port }to form a coal blend
that meets the specifications of a customer. Once a vessel arrives
at a berth at the terminal, the stockpiles with coal for the vessel
are reclaimed and loaded onto the vessel. The vessel then transports
the coal to its destination. The coordination of the logistics in
the Hunter Valley is challenging as it is a complex system involving
14 producers operating 35 coal mines, over 30 coal load points, 2 rail
track owners, 4 above rail operators, 3 coal loading terminals with
a total of 9 berths, and 9 vessel operators. Approximately 2000
vessels are loaded at the terminals in the Port of Newcastle each
year. For more information on the HVCC see the overview presentation
of the Hunter Valley Coal Chain Coordinator (HVCCC), the
organization responsible for planning the coal logistics in the
Hunter Valley \cite{HVCCC}.

\llong{We focus on the management of a stockyard at a coal loading
terminal\GLEB{, although we also propose model extensions looking at the processes outside the terminal}. An important characteristic of a coal loading terminal is
whether it operates as a cargo assembly terminal or as a dedicated
stockpiling terminal. When a terminal operates as a cargo assembly
terminal, it operates in a `pull-based' manner, where the coal
blends assembled and stockpiled are based on the demands of the
arriving ships. When a terminal operates as dedicated stockpiling
terminal, it operates in a `push-based' manner, where a small number
of coal blends are built in dedicated stockpiles and only these coal
blends can be requested by arriving vessels.

We focus on cargo assembly terminals as they are more difficult to manage
due to the large variety of coal blends that needs to be accommodated.
Depending on the size and the blend of a cargo, the assembly may take
anywhere from three to seven days. This is due, in part, to the fact that
mines can be located hundreds of miles away from the port and getting a
trainload of coal to the port takes a considerable amount of time. Once the
assembly of a stockpile has started, it is rare that the location of the
stockpile in the stockyard is changed; relocating a stockpile is
time-consuming and requires resources that can be used to assemble or
reclaim other stockpiles. Thus, deciding where to locate a stockpile and
when to start its assembly is critical for the efficiency of the
system. Ideally, the assembly of the stockpiles for a vessel completes at
the time the vessel arrives at a berth (i.e., `just-in-time' assembly) and
the reclaiming of the stockpiles commences immediately. Unfortunately, this
does not always happen due to the limited capacities of the resources in the
system, e.g., stockyard space, stackers, and reclaimers, and the complexity
of the stockyard planning problem.
}

\short{We focus on the management of a stockyard at one of the coal
loading terminals. It acts as a \emph{cargo assembly terminal} where
the coal blends assembled and stockpiled are based on the demands of
the arriving ships.} Our cargo assembly planning approach aims to
minimize the delay of vessels, where the delay of a vessel is
defined as the difference between the vessel's departure time and
its earliest possible departure time, that is, the departure time in
a system with infinite capacity. Minimizing the delay of vessels is
used as a proxy for maximizing the throughput, i.e., the maximum
number of tons of coal that can be handled per year, which is of
crucial importance as the demand for coal is expected to grow
substantially over the next few years. We investigate the value of
information given by the \emph{visibility horizon} --- the number of
future vessels whose arrival time and stockpile demands are known in
advance.


\llong{ Despite their economic importance, there is very little
literature on coal, and more generally mineral, supply chains.
Boland and Savelsbergh \cite{BSHVCC05} discuss a variety of planning
problems encountered in the HVCC. Singh et al.\ \cite{Singh12}
discuss expansion planning for the HVCC. Thomas et al.\
\cite{Thomas13} explore integrated planning and scheduling of coal
supply chains.
Gulzynsky et al.\ \cite{Gulcz12} present stockyard management
technology which combines greedy construction, enumeration, and
integer programming. The closest work to what we describe here is
that of \anonymize{Savelsbergh and Smith }\cite{SavSm13}, who propose a greedy
algorithm with partial lookahead and truncated tree search using
geometric properties of the space-time layouts, for the same
stockyard planning problem.  We propose a Constraint Programming (CP)-based approach and directly compare it with and extend
theirs. \llong{The core of our paper, namely the basic CP model and comparison with \cite{SavSm13}, was presented in the CPAIOR'14 paper \cite{BBSS14small}.}

There are obvious links between stockyard management and
2-dimensional bin or strip packing problems (see Lodi et al.\
\cite{LMV02} for a recent survey). The main difference is that in
stockyard management, the length of a rectangular item in one of its
dimensions (the time dimension) is not known in advance, but depends
on other decisions. Bay et al.\ \cite{Bay10} consider a shipping
yard planning problem where large 2-dimensional blocks need to be
handled. They solve this as a 3-dimensional bin packing problem in
which one dimension is time. }


The solving technology we apply is Constraint Programming using
Lazy Clause Generation (LCG)~\cite{lazyj}. Constraint programming
has been highly successful in tackling complex packing and
scheduling problems \cite{RCSPmax,Carpet11}. \llong{ CP allows
modeling optimization problems at higher level and the flexible
addition of new constraints without changing the base model; it also
allows the modeller to make use of their knowledge about where good
solutions are likely to be found by programming the search strategy. Lazy Clause
Generation \cite{LCG} is a new constraint programming solving
technology that allows the solver to ``learn'' from failures in the
search process. It often exponentially improves upon standard CP
technology. It is the state-of-the-art on many well-studied
scheduling problems and some packing problems.} Cargo assembly is a
combined scheduling and packing problem. The specific problem is
first described by \anonymRepl{Savelsbergh and Smith~}{[authors of]~}\cite{SavSm13}. They
propose a greedy heuristic for solving the problem and investigate
some options concerning various characteristics of the problem. We
present a Constraint Programming model implemented in the MiniZinc
language  \cite{MZN}. To solve the model efficiently, we develop
iterative solving methods: greedy methods to obtain initial
solutions and large neighbourhood search methods \cite{LNSPis} to
improve them. \llong{The CP technology facilitates investigation of
new practically relevant constraints and their impacts on the
solutions. \GLEB{In particular, we extend the model's scope by a port channel constraint which is modeled as a series of global $\disjunctive$ constraints}.  Finally, we propose an iterative adaptive scheme which
uses the greedy heuristic from \cite{SavSm13} as a core component and compare
both approaches. }





\llong{ The remainder of the paper is organized as follows.
Section~\ref{secModelBasic} describes the problem in detail and
reviews the new constraints we investigate with CP technology.
Section~\ref{secCP} describes the Constraint Programming model and
iterative methods based on it. Section~\ref{secAG} presents an
adaptive greedy scheme for the heuristic found in the literature.
Experiments are discused in Section~\ref{sec:exp} and finally in
Section~\ref{sec:conc} we conclude. }

\newcommand{\pA}{\textsf{A}}
\newcommand{\pB}{\textsf{B}}
\newcommand{\pC}{\textsf{C}}
\newcommand{\pD}{\textsf{D}}
\newcommand{\EH}{\textsf{EH}}
\newcommand{\LNS}{\textsf{LNS}}
\newcommand{\VH}{\textsf{VH}}

\section{Cargo assembly planning\llong{: overview}}\label{secModelBasic}

The starting point for this work is the model developed
in~\cite{SavSm13} for stockyard planning. \llong{The model is
described in the next subsection and summarised in
Subsection~\ref{secBAsicDescr}. We consider a number of constraints
beyond those in~\cite{SavSm13}, which are discussed in
Subsection~\ref{secFurtherConstrOverview}.


\subsection{Cargo assembly environment}  \label{secCAEnv}
}

\begin{figure}[tb]
\unitlength0.0004\textwidth\centering
\begin{picture}(2200,450)
\put(58,350){\color[RGB]{200,200,200}\rule{2142\unitlength}{70\unitlength}}
\put(1100,385){\makebox(0,0){Pad \pA}}
\put(295,240){\color[RGB]{200,200,200}\rule{1905\unitlength}{70\unitlength}}
\put(1100,275){\makebox(0,0){Pad \pB}}
\put(26,130){\color[RGB]{200,200,200}\rule{2174\unitlength}{70\unitlength}}
\put(1100,165){\makebox(0,0){Pad \pC}}
\put(44,20){\color[RGB]{200,200,200}\rule{2156\unitlength}{70\unitlength}}
\put(1100,55){\makebox(0,0){Pad \pD}}
\newsavebox{\svbox}
\savebox{\svbox}(0,0)[bl]{
\put(0,0){\rule{560\unitlength}{4\unitlength}}
\put(740,0){\rule{630\unitlength}{4\unitlength}}
\put(1540,0){\rule{610\unitlength}{4\unitlength}}
}
\put(50,440){\usebox{\svbox}}
\put(700,440){\makebox(0,0){\footnotesize\textsl{S1}}}
\put(1500,440){\makebox(0,0){\footnotesize\textsl{S2}}}
\put(50,330){\usebox{\svbox}}
\put(700,330){\makebox(0,0){\footnotesize\textsl{R1}}}
\put(1500,330){\makebox(0,0){\footnotesize\textsl{R2}}}
\put(50,220){\usebox{\svbox}}
\put(700,220){\makebox(0,0){\footnotesize\textsl{S3}}}
\put(1500,220){\makebox(0,0){\footnotesize\textsl{S4}}}
\put(50,110){\usebox{\svbox}}
\put(700,110){\makebox(0,0){\footnotesize\textsl{R3}}}
\put(1500,110){\makebox(0,0){\footnotesize\textsl{R4}}}
\put(50,0){\usebox{\svbox}}
\put(700,0){\makebox(0,0){\footnotesize\textsl{S5}}}
\put(1500,0){\makebox(0,0){\footnotesize\textsl{S6}}}
\end{picture}
\caption{A scheme of the stockyard with 4 pads, 6 stackers, and 4 reclaimers}\label{figKCT}
\end{figure}
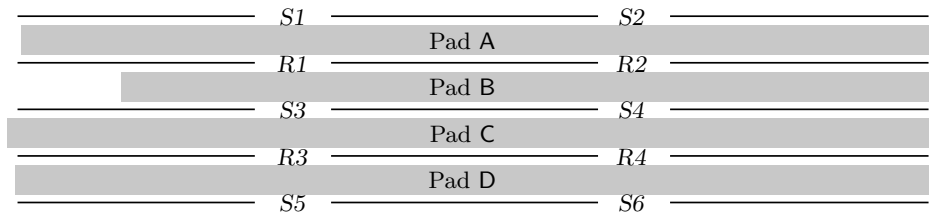

The stockyard studied has four pads, \pA, \pB, \pC,
and \pD, on which cargoes are assembled. Coal arrives at the
terminal by train. Upon arrival at the terminal, a train dumps its
contents at one of three dump stations. The coal is then transported
on a conveyor to one of the pads where it is added to a stockpile by
a stacker. There are six stackers, two that serve pad \pA, two that
serve both pads \pB{} and \pC, and two that serve pad \pD. A single
stockpile is built from several train loads over several days. After
a stockpile is completely built, it dwells on its pad for some time
(perhaps several days) until the vessel onto which it is to be
loaded is available at one of the berths. A stockpile is reclaimed
using a bucket-wheel reclaimer and the coal is transferred to the
berth on a conveyor. The coal is then loaded onto the vessel by a
shiploader. There are four reclaimers, two that serve both pads
\pA{} and \pB{}, and two that serve both pads \pC{} and \pD. Both
stackers and reclaimers travel on rails at the side of a pad.
Stackers and reclaimers that serve the same pads cannot pass each
other. A scheme of the stockyard is given in Figure~\ref{figKCT}.

\short{ The cargo assembly planning process involves the following
steps. An incoming vessel defines a set of cargoes (different blends
of coal) to be assembled and an \emph{estimated time of arrival
(ETA)}. The cargoes are assembled in the stockyard as different
stockpiles. The vessel cannot arrive at berth earlier than its ETA.
Once at a berth, and once all its cargoes have been assembled, the
reclaiming of the stockpiles (the loading of the vessel) begins. The
stockpiles are reclaimed onto the vessel in a specified order to
maintain physical balancing constraints. The goal of the planning
process is to maximize the throughput without causing unacceptable
delays for the vessels.

When assigning each cargo of a vessel to a location in the stockyard
we need to schedule the stacking and reclaiming of the stockpile
taking into account limited stockyard space, stacking rates,
reclaiming rates, and reclaimer movements. We model stacking and
reclaiming at different levels of granularity. All reclaimer
activities, e.g., the reclaimer movements along its rail track and
the reclaiming of a stockpile, are modelled in time units of one
minute. Stacking is modelled only at a coarse level of detail in 12
hour periods.

We assume that the time to build a stockpile is derived from the
locations of the mines that contribute coal to the blend (the
distance of the mines from the port). We allocate 3, 5, or 7 days to
stacking of different stockpiles depending on the blend. We assume
that the tonnage of the stockpile is stacked evenly over the
stacking period. Since the trains that transport coal from the mines
to the terminal are scheduled closer to the day of operations, this
is not unreasonable. We assume that all stockpiles for a vessel are
assembled on the same pad, since that leads to better results
(already observed in \cite{SavSm13}). In practice, however, there is
no such restriction.

For each stockpile we need to decide a location, a stacking start
time, a reclaiming start time, and which reclaimer will be used.
Note that reclaiming does not have to start as soon as stacking has
finished; the time between the completion of stacking and the start
of reclaiming is known as \emph{dwell time}. Stockpiles cannot
overlap in time and space, reclaimers can only be used on pads they
serve, and reclaimers cannot cross each other on the shared track.
The waiting time between the reclaiming of two consecutive
stockpiles of one vessel is limited by the \emph{continuous reclaim
time limit}.
The reclaiming of a stockpile, a so-called \emph{reclaim job},
cannot be interrupted.

A cargo assembly plan can conveniently be represented using
\emph{space-time diagrams}; one space-time diagram for each of the
pads in the stockyard. A space-time diagram for a pad shows for any
point in time which parts of the pad are occupied by stockpiles (and
thus also which parts of the pad are not occupied by stockpiles and
are available for the placement of additional stockpiles) and the
locations of the reclaimers serving that pad. Every pad is
rectangular; however its width is much smaller than its length and
each stockpile is spread across the entire width. Thus, we model
pads as one-dimensional entities. The location of a stockpile can be
characterized by the position of its lowest end called its
\emph{height}.
A stockpile occupies space on the pad for a certain amount of time. This time can be divided
into three distinct parts: the \emph{stacking part}, i.e., the time during which the stockpile is being
built; the \emph{dwell part}, i.e., the time between the end of stacking and the start of reclaiming;
and a \emph{reclaiming part}, i.e., the time during which the stockpile is reclaimed and loaded on a
waiting vessel at a berth. Thus, each stockpile can be represented in a space-time diagram
by a three-part rectangle as shown in Figure \ref{figSTDiag}.
\begin{figure}[tb]\centering
\includegraphics[page=12,width=0.9\textwidth]{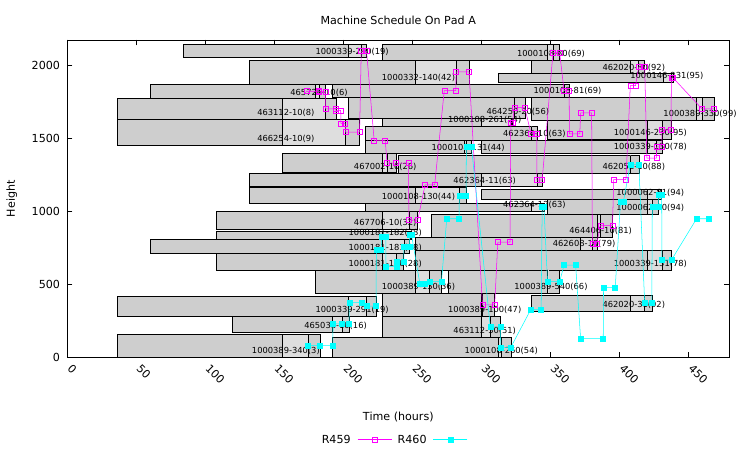}
\caption{A space-time diagram of pad A showing also reclaimer movements. Reclaimer R459 has to be after R460 on the pad. Both reclaimers also have jobs on pad B.}\label{figSTDiag}
\end{figure}} \llong{ A brief overview of the events
driving the cargo assembly planning process is presented next. An
incoming vessel alerts the coal chain managers of its pending
arrival at the port. This announcement is referred to as the
vessel's nomination. Upon nomination, a vessel provides its
\emph{estimated time of arrival (ETA)} and a specification of the
cargoes to be assembled to the coal chain managers. As coal is a
blended product, the specification includes for each cargo a recipe
indicating from which mines coal needs to be sourced and in what
quantities (this is beyond our current scope however). At this time, the assembly of the cargoes (stockpiles)
for the vessel can commence. A vessel cannot arrive at a berth prior
to its ETA, and often a vessel has to wait until after its ETA for a
berth to become available. Once at a berth, and once all its cargoes
have been assembled, the reclaiming of the stockpiles (the loading
of the vessel) can begin. A vessel must be loaded in a way that
maintains its physical balance in the water. As a consequence, for
vessels with multiple cargoes, there is a predetermined sequence in
which its cargoes must be reclaimed. The goal of the planning
process is to maximize the throughput without causing unacceptable
delays for the vessels.

For a given set of vessels arriving at the terminal, the goal is thus to assign each cargo of
a vessel to a location in the stockyard, schedule the assembly of these cargoes, and schedule
the reclaiming of these cargoes, so as to minimize the average delay of the vessels, where the
delay of a vessel is defined to be the difference between the departure time of the vessel (or
equivalently the time that the last cargo of the vessel has been reclaimed) and the earliest
time the vessel could depart under ideal circumstances, i.e., the departure time if we assume
the vessel arrives at its ETA and its stockpiles can be reclaimed immediately upon its arrival.

When assigning each cargo of a vessel to a location in the stockyard, scheduling the
assembly of these cargoes, and scheduling the reclaiming of these cargoes, we account for
limited stockyard space, stacking rates, reclaiming rates, and reclaimer movements, but we
assume that other parts of the system have infinite capacity and thus do not restrict our
decisions. The reason for doing so is that our industry partner believes that the stockyard
is, or will soon be, the bottleneck in the system and that the other parts of the system can
be managed in such a way that they adjust to what is needed to operate the stockyard as
effectively as possible. Consequently, we assume that coal is available for stacking any time,
and that a vessel is at berth and can be loaded whenever we want to reclaim a stockpile.

We model stacking and reclaiming at different levels of granularity. Since reclaimers are
most likely to be the constraining entities in the system, we model reclaimer activities at a
fine level of detail. That is, all reclaimer activities, e.g., the reclaimer movements along its
rail track and the reclaiming of a stockpile, are modelled in time units of one minute. On the other
hand, we model stacking at a coarser level of detail. First, we assume that the time it
takes to build a stockpile can be derived from its constituent components. As the constituent
components define the mines from which the coal will be sourced, this is not unreasonable.
The distance from the terminal to the mine determines the cycle time of the trains used
to transport the coal and thus stockpiles with coal sourced from mines that are far away
from the terminal tend to take longer to build. Second, we consider stacking capacity at an
aggregate daily level for each pair of stackers serving the same pads. As each pair of stackers
is linked by a conveyor system to a particular dump station at the terminal, often referred
to a stacking stream, this is not unreasonable.

When deciding a stockpile location, a stockpile stacking start time, and a stockpile reclaiming start time, a number of constraints have to be taken into account: at any point in
time no two stockpiles can occupy the same space on a pad, reclaimers cannot be assigned
to two stockpiles at the same time, reclaimers can only be assigned to stockpiles on pads
that they serve, reclaimers serving the same pads cannot pass each other, the stockpiles of a
vessel have to be reclaimed in a specified reclaim order and the time between the reclaiming
of consecutive stockpiles of a vessel can be no more than a prespecified limit, the so-called
\emph{continuous reclaim time limit}, and, finally, the reclaiming of the first stockpile of a vessel
cannot start before all stockpiles of that vessel have been stacked. We assume that the rate
of reclaiming is the same for all reclaimers, and thus that the time it takes to reclaim a given
stockpile can be calculated from its size, which is a reasonable model of the real situation.
\didpjs

A cargo assembly plan can conveniently be represented using
\emph{space-time diagrams}; one space-time diagram for each of the
pads in the stockyard. A space-time diagram for a pad shows for any
point in time which parts of the pad are occupied by stockpiles (and
thus also which parts of the pad are not occupied by stockpiles and
are available for the placement of additional stockpiles) and the
locations of the reclaimers serving that pad. Every pad is
rectangular; however its width is much smaller than its length and
each stockpile is spread across the entire width. Thus, we model
pads as one-dimensional entities. The location of a stockpile can be
characterized by the position of its lowest end called its
\emph{height}.
A stockpile occupies space on the pad for a certain amount of time. This time can be divided
into three distinct parts: the \emph{stacking part}, i.e., the time during which the stockpile is being
built; the \emph{dwell part}, i.e., the time between the end of stacking and the start of reclaiming;
and a \emph{reclaiming part}, i.e., the time during which the stockpile is reclaimed and loaded on a
waiting vessel at a berth. We assume that a stockpile's reclaim process, \emph{reclaim job},\NEW{} cannot be interrupted. Thus, each stockpile can be represented in a space-time diagram
by a three-part rectangle as shown in Figure \ref{figSTDiag}.
\begin{figure}[tb]\centering
\includegraphics[page=12,width=0.9\textwidth]{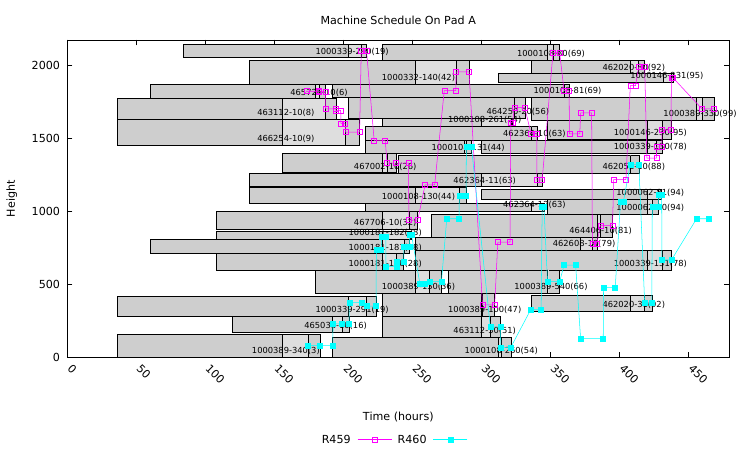}
\caption{A space-time diagram of pad A showing also reclaimer movements. Reclaimer R459 has to be after R460 on the pad. Both reclaimers also have jobs on pad B.}\label{figSTDiag}
\end{figure}

\subsection{The basic model}  \label{secBAsicDescr}

Based on the above description of the problem, we summarize the important
features of the model studied in \cite{SavSm13}, which is taken as the basic model here.
\begin{itemize}\parskip0em\parsep0em\itemsep0em
\item 4 pads, arranged in parallel, 4 reclaimers (2 between each pair of pads).
\item 3 stacker streams, (one each for the outer pads, one for the two inner), each with capacity 288,000 t/day and overall capacity (maximum \emph{daily inbound throughput, DIT}) 537,600 t/day.
\item Stacking duration 3/5/7 days for each stockpile, which is reflective of what happens in real life.
\item The stacking volume of each stockpile is evenly divided over the stacking period.
\item Stacking of a stockpile can start at most 10 days before the ETA of
  its vessel.
\item Reclaimer's travel speed: 30 meters per minute.
\item Reclaim the stockpiles of each vessel in given order, one stockpile at a time.
\item Reclaim job of any stockpile is non-interruptible.\NEW
\item Continuous reclaim limit (i.e., maximal time between reclaiming) for two consecutive stockpiles of a vessel: 5 hours. \didpjs
\item Same pad is used for all stockpiles of a single
  vessel. \didpjs
  While not the industrial practice, this rule produced better solutions both in \cite{SavSm13} and in our experiments, so we include it into the basic model and present its impact in the numerical results.
\end{itemize}
The goal is to minimize the average delay of all vessels.
 \par In \cite{SavSm13} the authors do not consider the delay of the first and last 8 vessels. The reason is that these vessels have more placement freedom and thus easy to schedule. In our experiments, exclusion of the delays in these `warm-up' and `cool-down' vessels from the objective function made them unreasonably high.\NEW{} That is why we optimize over all vessels, however we report some results taking `warm-up' and 'cool-down' into account. Moreover, we investigate the effect of a limited visibility horizon, i.e., when the number of the next arriving ships known is limited.

As observed in \cite{SavSm13}, minimizing dwell (stockpile's waiting time) after scheduling each stockpile or vessel is generally bad because this reduces resources for later vessels. Thus, we do not consider dwell in the basic model.


\subsection{\GLEB{Model extensions}} \label{secFurtherConstrOverview}
There are many further constraints which are relevant in practice and it is important to investigate their effect on the solutions. We consider the following optional constraints\GLEB{, some of which are a first step toward a system view of the transport chain}:
\begin{itemize}\parskip0em\parsep0em\itemsep0em
\item Different pads allowed for allocation of stockpiles of one and the same vessel.\NEW
\item Build all stockpiles of a vessel (stack them) before reclaiming any of them.
\item Maximal number of simultaneously berthed ships (4).
\item Maximal number of reclaimers working at any time (3). This constraint models the fact that only 3 ship loaders are available at the terminal.\NEW
\item Tidal constraints: big ships are tidally constrained, that is, they
can leave the port only during high tides, 
which for simplicity we model as the periods 
11:15--12:45am and 11:15--12:45pm (more accurate modeling would certainly
not be too difficult to add). In practice, a ship is considered tidally constrained if its weight is equal to or above 100,000 tonnes.

\item Channel constraints: time between any two departures is at least 20 minutes; the same for arrivals.\NEW{} After any departure, the earliest possible arrival is 140 minutes later. The reason is the port entry channel which can let only one vessel in each direction.
\item Flexible stacking volumes for each day for a stockpile.
\item $\pm$1 day stacking duration for each stockpile.
\item Dwell post-processing. 
\end{itemize}
The next section presents a Constraint Programming model of the above problem.


%
%
}

\llong{
\section{The Constraint Programming-based approach}\label{secCP}


We implemented the model in the MiniZinc language \cite{MZN} which is
accepted by many solvers.
As reasoned in Section~\ref{secCAEnv}, we use a
fine and coarse time discretization for modeling reclaiming and stacking, respectively.
The
non-overlapping of stockpiles in space and time is a three-dimensional
\texttt{diffn} constraint \cite{BeldiceanuGlobal94} where the 
pad number can
be seen as a third dimension.
In the solver we use the \texttt{diffn} constraint is implemented by
decomposition.
}

\llong{
The model extensively uses the constraint
\be
\cumulative(t,d,r,R)
\ee
restricting the available renewable resource capacity $R$ for jobs with
start times, durations, and resource demands given by vectors $t$, $d$, and
$r$, respectively \cite{MZN}.
If $d$, $r$, and $R$ are constants, the particular
solver we use applies a \emph{global} version of $\cumulative$,
otherwise a decomposed version \cite{CumulExpl11} is used
(e.g. for constraint \eqref{eqMaxBerth}).
A redundant $\cumulative$ constraint \eqref{eqCumPad} models pad usage, which is known to be
important to reduce search space in packing problems \cite{Carpet11}.
Cumulative constraints are also used to model stacking capacity and,
in the extended model (Section~\ref{secFurtherConstr}), reclaimer usage. 

We also make use of the $\disjunctive$ global constraint which
is a specialized form of $\cumulative$.
\be
\disjunctive(t,d) = \cumulative(t,d,\mathbf{1}, 1)
\ee
where $\mathbf{1}$ is a vector of 1s,
which ensures that no two tasks overlap in time.
Constraint programming solvers can more efficiently handle 
the specialized $\disjunctive$ constraint.



We present the basic model, extended model, solver search strategy, construction of initial solutions, large neighbourhood search, and an approach for limited visibility horizons.
}

\subsection{The basic Constraint Programming model}\label{secModelBasicCP}

\llong{ This subsection models the basic problem subset described in
Section~\ref{secBAsicDescr}. The presented structure of the
constraints corresponds to their implementation in the MiniZinc
language, however mathematical notation is used where possible to
improve readability. } \short{ We present the model of the cargo
assembly problem below; the structure corresponds directly to the
implementation in MiniZinc~\cite{minizinc}. The unit for time
parameters is minutes, and for space parameters is meters.}
\llong{In \cite{SavSm13} continuous input data is used.\NEW{} But Constraint Programming solvers concentrate on integer
variables, so in the implementation the times are rounded to
minutes; distances to whole meters; and stacking tonnages to
integers. We round up in order to be conservative.} 
In addition, stacking start times are restricted to be multiples of 12 hours.

\newlength{\defssLMD}
\newlength{\defssLW}
\defssLMD8.8em
\defssLW10.2em
\newenvironment{defss}
{
    \begin{list}{}{\parsep0em\itemsep0em\parskip0em%
    \addtolength{\leftmargin}{\defssLMD}
    \labelwidth\defssLW
    }
}{
    \end{list}
}
\newcommand{\defssitem}[1]{\item[{\makebox[\defssLW][l]{{#1}\hfill---\ }}]}
\newcommand{\defsalone}[1]{\item[{\makebox[\defssLW][l]{{#1}}}]}

\begin{defss}
\defsalone{\textbf{Parameter sets}}
\defssitem{$S$} set of stockpiles of all vessels, ordered by vessels' ETAs and reclaim sequence of each vessel's stockpiles
\defssitem{$V$} set of vessels, ordered by ETAs
\end{defss}


\begin{defss}
\defsalone{\textbf{Parameters}}
\defssitem{$v_s$} vessel for stockpile $s \in S$
\defssitem{$\veta_v$} estimated time of arrival of vessel $v\in V$, minutes
\item[$d^S_s\in\{4320,7200,10080\}$ \quad --- \ \;stacking duration of stockpile $s\in S$, minutes]
\defssitem{$d^R_s$} reclaiming duration of stockpile $s\in S$, minutes
\defssitem{$l_s$} length of stockpile $s\in S$, meters
\item[$(H_1,\dots,H_4) = (2142, 1905, 2174, 2156)$ \qquad --- \ \;pad lengths, meters]
\defssitem{$\speed^R\ =30$} travel speed of a reclaimer, meters / minute
\defssitem{$\tonn^\text{daily}_s$} daily stacking tonnage of stockpile $s\in S$, tonnes
\defssitem{$\tonn^\text{DIT}=\text{537,600}$} daily inbound throughput (total daily stacking capacity), tonnes
\defssitem{$\tonn^\text{SS \ }_k=\text{288,000}$} daily capacity of stacker stream $k\in\{1,2,3\}$, tonnes
\end{defss}

\defssLMD11.7em
\defssLW13.5em
\begin{defss}
\defsalone{\makebox[5.5em]{\textbf{Decisions}}}
\defssitem{\makebox[2em]{$p_v$}$\in\{1,\dots,4\}$} pad on which the stockpiles of vessel $v\in V$ are assembled
\defssitem{\makebox[2em]{$h_s$}$\in\{0,\dots,H_{p_{v_s}}-l_s\}$} position of stockpile $s\in S$ (of its `closest to pad start' boundary) on the pad
\defssitem{\makebox[2em]{$t^S_s$}$\in\{0,720,\dots\}$} stacking start time of stockpile $s\in S$
\defssitem{\makebox[2em]{$r_s$}$\in\{1,\dots,4\}$} reclaimer used to reclaim stockpile $s\in S$
\defssitem{\makebox[2em]{$t^R_s$}$\in\{\veta_{v_s},\veta_{v_s}+1,\dots\}$} reclaiming start time of stockpile $s\in S$
\end{defss}

\subsubsection*{Constraints.}

\llong{\begin{subequations}}

\noindent Reclaiming of a stockpile cannot start before its
vessel's ETA: \be\nonumber t^R_s \ge \veta_{v_s},\qquad \forall s
\in S \ee \llong{\\}
Stacking of a stockpile starts no more than 10 days before its
vessel's ETA: \be\nonumber t^S_s \ge \veta_{v_s}-14400,\qquad
\forall s \in S \ee \llong{\\}
Stacking of a stockpile has to complete before reclaiming can start:
\be\nonumber t^S_s + d^S_s \le t^R_s, \qquad \forall s \in S \ee \llong{\\}
The reclaim order of the stockpiles of a vessel has to be respected:
\be\nonumber t^R_s + d^R_s \le t^R_{s+1},\qquad \forall s \in S
\text{~where~} \ v_s=v_{s+1} \ee \llong{\\}
The continuous reclaim time limit of 5 hours has to be respected:
\be\nonumber t^R_{s+1} - 300 \le t^R_s + d^R_s,\qquad \forall s \in
S \text{~where~} \ v_s=v_{s+1} \ee \llong{\\}
A stockpile has to fit on the pad it is assigned to:
\be\nonumber 0\le h_s\le H_{p_{v_s}}-l_s, \qquad \forall s \in S \ee \llong{\\}
Stockpiles cannot overlap in space and time:
\begin{multline}
p_{v_s}\ne p_{v_t} \vee h_s+l_s\le h_t\vee h_t+l_t\le h_s
    \vee t^R_s+d^R_s\le t^S_t\vee t^R_t+d^R_t\le t^S_s, \\ \forall s<t\in S \nonumber
\end{multline}
Reclaimers can only reclaim stockpiles from the pads they serve:
\be\nonumber p_{v_s}\le2 \Leftrightarrow r_s\le2, \qquad \forall s \in S \ee \llong{\\}
If two stockpiles $s<t$ are reclaimed by the same reclaimer,
then the time between the end of reclaiming the first and the start
of reclaiming the second should be enough for the reclaimer to move
from the middle of the first to the middle of the second:
\be\label{eqReclMov}\short{\nonumber}
r_s\ne r_t
 \vee \max\bigl\{(t^R_t-t^R_s-d^R_s),(t^R_s-t^R_t-d^R_t)\bigr\}\speed^R \ge \Bigl|h_s+\frac{l_s}2-h_t-\frac{l_t}2\Bigr| \short{\nonumber} %
\ee
To avoid clashing, at any point in time, the position of Reclaimer~2
should be before the position of Reclaimer~1 and the position of
Reclaimer~4 should be before the position of Reclaimer~3. An example
of the position of Reclaimers 1 and 2 in space and time is given in
Figure~\ref{figClAv} (see also Figure~\ref{figSTDiag}).
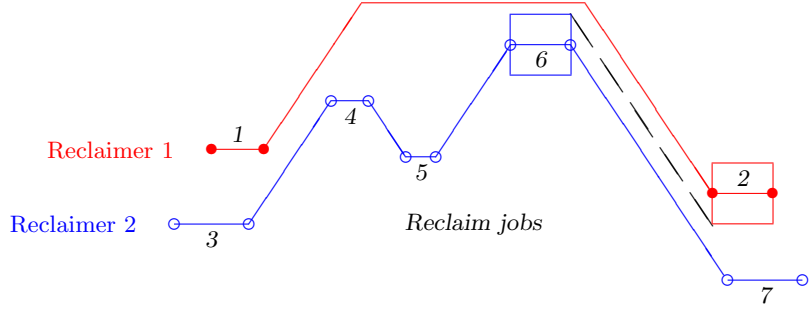
\begin{figure}[tb]
\unitlength0.007\textwidth\centering
\begin{picture}(100,40)
\put(60,10){\makebox(0,0)[B]{\sl Reclaim jobs}}
\multiput(73,38.2)(4,-6){5}{\line(2,-3){3}}
\color{red}
\put(20,20){\makebox(0,0)[r]{Reclaimer 1}}
\put(25,20){\put(0,0){\circle*{1.5}}\put(0,0){\line(1,0){7}}\put(7,0){\circle*{1.5}}\put(0,1){\makebox(7,0)[b]{\color{black}\sl1}}}
\put(32,20){\line(2,3){13}}
\put(45,39.6){\line(1,0){30}}
\put(75,39.5){\line(2,-3){17}}
\put(92,14.1){\put(0,0){\circle*{1.5}}\put(0,0){\line(1,0){8}}\put(8,0){\circle*{1.5}}\put(0,1){\makebox(8,0)[b]{\color{black}\sl2}}}
\put(92,10.1){\framebox(8,8){}}
\color{blue}
\put(15,10){\makebox(0,0)[r]{Reclaimer 2}}
\put(20,10){\put(0,0){\circle{1.5}}\put(0,0){\line(1,0){10}}\put(10,0){\circle{1.5}}\put(0,-1){\makebox(10,0)[t]{\color{black}\sl3}}}
\put(30,10){\line(2,3){11}}
\put(41,26.5){\put(0,0){\circle{1.5}}\put(0,0){\line(1,0){5}}\put(5,0){\circle{1.5}}\put(0,-1){\makebox(5,0)[t]{\color{black}\sl4}}}
\put(46,26.5){\line(2,-3){5}}
\put(51,19){\put(0,0){\circle{1.5}}\put(0,0){\line(1,0){4}}\put(4,0){\circle{1.5}}\put(0,-1){\makebox(4,0)[t]{\color{black}\sl5}}}
\put(55,19){\line(2,3){10}}
\put(65,34){\put(0,0){\circle{1.5}}\put(0,0){\line(1,0){8}}\put(8,0){\circle{1.5}}\put(0,-1){\makebox(8,0)[t]{\color{black}\sl6}}}
\put(65,30){\framebox(8,8){}}
\put(73,34){\line(2,-3){20.5}}
\put(94,2.5){\put(0,0){\circle{1.5}}\put(0,0){\line(1,0){10}}\put(10,0){\circle{1.5}}\put(0,-1){\makebox(10,0)[t]{\color{black}\sl7}}}
\end{picture}
\caption{A schematic example of space (vertical)-time (horizontal) location of Reclaimers 1 and 2 with some reclaim jobs. Reclaimer~2 has to stay spatially before Reclaimer~1.}\label{figClAv}
\end{figure}
\llong{In Figure~\ref{figClAv}, because}\short{Because} Job~3 is
spatially before Job~1, there is no concern for a clash. However,
since Job~6 is spatially after Job~2, we have to ensure that there
is enough time for the reclaimers to get out of each other's way.
The slope of the dashed line corresponds to the reclaimer's travel
speed ($\speed^R$), so we see that the time between the end of Job 6
and the start of Job 2 has to be at least $(h_6 + l_6 -
h_2)/\speed^R$.

We model clash avoidance by a disjunction: for any two stockpiles $s
\ne t$, one of the following conditions must be met: either $(r_s
\ge 3 \wedge r_t \le 2)$, in which case $r_s$ and $r_t$ serve
different pads; or $r_s < r_t$, in which case $r_s$ does not have to
be before $r_t$; or $h_s + l_s \le h_t$, in which case stockpile $s$
is before stockpile $t$; or, finally, enough time between the reclaim
jobs exists for the reclaimers to get out of each other's way:
\begin{multline}\llong{\label{eqReclClash}}
\short{\nonumber}
 \max\bigl\{(t^R_t-t^R_s-d^R_s),(t^R_s-t^R_t-d^R_t)\bigr\}\speed^R \ge h_s+l_s-h_t
 \\\short{\nonumber}
\vee r_s < r_t \vee (r_s\ge3 \wedge r_t\le2) \vee h_s+l_s\le h_t, \quad \forall s\ne t\in S
\end{multline}
Redundant cumulatives on pad space usage improved efficiency. They
require derived variables $l^p_s$ giving the `pad length of
stockpile $s$ on pad $p$':
\begin{gather}\nonumber
l^p_s=\begin{cases}l_s, &\text{if } p_{v_s}=p,
\\0, &\text{otherwise}, \end{cases}
 \qquad \forall s\in S, \ p\in\{1,\dots,4\}
\\\label{eqCumPad}\short{\nonumber}
\cumulative(t^S, t^R+d^R-t^S, l^p, H_p), \qquad p\in\{1,\dots,4\}
\end{gather}
The stacking capacity is constrained day-wise. If a stockpile is
stacked on day $d$ and the stacking is not finished before the end
of $d$, the full daily tonnage of that stockpile is accounted for
using derived variables \be\nonumber t^{S1}=\lfloor
t^S/1440\rfloor,\quad d^{S1}=\lfloor d^S/1440\rfloor \ee The daily
stacking capacity cannot be exceeded: \be\nonumber
\cumulative(t^{S1}, d^{S1}, \tonn^\text{daily}, \tonn^\text{DIT})
\ee The capacity of stacker stream $k$ (a set of two stackers
serving the same pads) is constrained similar to pad space usage:
\begin{gather}\label{eqStStream}\nonumber
\tonn^\text{daily}_{ks}=\begin{cases}\tonn^\text{daily}_s, &\text{if } (p_{v_s},k)\in\{(1,1),(2,2),(3,2),(4,3)\}
\\0, &\text{otherwise}, \end{cases}
 \qquad \forall s,k
\\\nonumber
\cumulative(t^{S1}, d^{S1}, \tonn^\text{daily}_k, \tonn^\text{SS}_k), \qquad k\in\{1,2,3\}
\end{gather}
\short{The maximum number of simultaneously berthed ships is 4. We
introduce derived variables for vessels' berth arrivals and use a
decomposed cumulative:
\begin{gather}\label{eqBerth}\nonumber
t^\text{Berth}_v=t^R_{s^\text{first}(v)}, \llong{\qquad}\short{\qquad s^\text{first}(v) = \min\{s | v_s=v\},\qquad\quad }\forall v\in V\\\nonumber
\card(\{u\in V\ |\ u\ne v, \ t^\text{Berth}_u\le t^\text{Berth}_v
    \wedge t^\text{Berth}_v < t^\text{Depart}_u\}) \le 3,\qquad \forall v\in V  \label{eqMaxBerth}
\end{gather}
}

\subsubsection*{Objective function.}

The objective is to minimize the sum of vessel delays. To define
vessel delays, we introduce the derived variables
$t^\text{Depart}_v$ for vessel departure times:
\begin{align} \llong{\label{eqDepEarl}}\short{\nonumber}
\depEarliest_v &= \veta_v+\sum_{s | v_s=v} d^R_s,\qquad &\forall v\in V
\\\llong{\label{eqTDepart}}\short{\nonumber}
t^\text{Depart}_v &= t^R_{s^\text{last}(v)} + d^R_{s^\text{last}(v)}, \llong{\qquad} \short{\quad}&s^\text{last}(v) = \max\{s | v_s=v\},\quad\forall v\in V
\\\llong{\label{eqDelay}}\short{\nonumber}
\delay_v &= t^\text{Depart}_v-\depEarliest_v,\qquad &\forall v\in V
\\\textstyle\label{eqObj}
\textbf{objective} &= \sum_v \delay_v
\end{align}
\llong{\end{subequations}}
\llong{
Except for the discretizations\short{ and the berth capacity constraint},
the above model corresponds to that used in
\cite{SavSm13}.}
\llong{
The differences are once more summarized in Section~\ref{secDiff}.



\subsection{\GLEB{Model extensions}}\label{secFurtherConstr}
Constraint programming models allow relatively
easy addition of further constraints and options to the model\GLEB{, either detailing stockyard operation or leading to a more global view of the system}.
Below we explain and define them. 

\begin{subequations}
\paragraph{Different pads for stockpiles of one vessel allowed:} instead of the variables $p_v$ denoting the pad used to allocate the stockpiles of vessel $v\in V$,\NEW{} implement $p_s$ for each stockpile $s\in S$.
\paragraph{Stack all before reclaim:} Build all stockpiles of a vessel (stack them) before reclaiming any of them.
\begin{multline}\label{eqSAbR}
t^S_s+d^S_s \le t^R_{s^\text{first}(v)}, \short{\\} \llong{\quad} \forall v_s=v:s\ne s^\text{first}(v)=\min\{s|v_s=v\}, \ \forall v\in V
\end{multline}
\paragraph{Berth capacity:} Maximal number of simultaneously berthed ships is 4.
We introduce derived variables for vessels' berth arrivals and use a decomposed cumulative:
\begin{gather}\label{eqBerth}
t^\text{Berth}_v=t^R_{s^\text{first}(v)}, \qquad\forall v\in V\\
\card(\{u\in V\ |\ u\ne v, \ t^\text{Berth}_u\le t^\text{Berth}_v
    \wedge t^\text{Berth}_v < t^\text{Depart}_u\}) \le 3=4-1,\qquad \forall v\in V  \label{eqMaxBerth}
\end{gather}

\paragraph{Ship loading capacity:} Maximal number of reclaimers working at any time is 3.
This can be easily imposed by a global cumulative:
\be\label{eqMaxReclaim}
\cumulative(t^R, d^R, \mathbf1, 3)
\ee

\paragraph{Tidal constraints:} Ships exceeding certain weight limit can leave the port only during the high tides 11:15--12:45am and 11:15--12:45pm. The berth is allocated until departure.

Modeling: the departure time variables $t^\text{Depart}_v$ \eqref{eqTDepart} for tidal vessels are disconnected from the end of reclaiming by changing the equalities \eqref{eqTDepart} to `greater-or-equal'. The domains of these variables are set as the tidal windows.\NEW{} 
  Moreover, the constant representing earliest possible departure
  \eqref{eqDepEarl} and used to compute actual delay \eqref{eqDelay} is
  increased to the next tidal window boundary if it is not already inside a tidal window. 

\paragraph{Channel constraint:} time between any two departures is at least
20 minutes; the same for arrivals.\NEW{} After any departure, the earliest
possible arrival is 140 minutes later. The channel constraint can be
illustrated in a space-time diagram representing vessels' positions in the
channel depending on time. To pass through the channel in any direction, a
vessel needs 60 minutes. Moreover, before a new vessel can enter, the last
exiting vessel needs at least 20 more minutes to clear the entrance area. We
can represent arrivals and departures by polygons in a space-time
diagram. The left-hand side of a polygon represents a vessel's movement. The
polygons have time spans 60 and 80 for arrival and departure, respectively,
and `thickness' 20 to ensure the distance between vessels going in one direction. This gives a polygon packing problem exemplified in Figure~\ref{figChannel}. We don't use this analogy for our modeling, however.
\begin{figure}[tb]
\unitlength0.0015\textwidth\centering
\begin{picture}(550,200)
\put(20,60){\rotatebox{90}{\makebox(0,0)[rb]{\parbox{6ex}{\raggedleft Ch.\\entry\\area}}}}
\thicklines
\put(30,30){\line(0,1){140}}
\put(20,65){\line(1,0){20}}
\put(20,92.5){\makebox(0,0)[rb]{\rotatebox{90}{Channel}}}
\put(30,30){\vector(1,0){530}}
\put(30,170){\vector(1,0){530}}
\thinlines
\newsavebox{\svarr}
\savebox{\svarr}(0,0)[bl]{
\put(0,0){\line(0,1){35}}
\put(0,0){\line(2,3){70}}
\put(70,105){\line(0,1){35}}
\put(0,35){\line(2,3){70}}
}
\newsavebox{\svdep}
\savebox{\svdep}(0,0)[bl]{
\put(0,140){\line(1,0){23.333}}
\put(0,140){\line(2,-3){93.333}}
\put(93.333,0){\line(0,1){35}}
\put(23.333,140){\line(2,-3){70}}
}
\put(80,30){\usebox{\svarr}}
\put(150,30){\usebox{\svdep}}
\put(173.333,30){\usebox{\svdep}}
\put(266.666,30){\usebox{\svarr}}
\put(290,30){\usebox{\svarr}}
\put(360,30){\usebox{\svdep}}
\put(453.333,30){\usebox{\svarr}}
\put(75,20){\rotatebox{-25}{\makebox(0,0)[l]{$t^A_1-60$}}}
\put(142,178){\rotatebox{45}{\makebox(0,0)[l]{$t^A_1=t^D_2$}}}
\put(171.333,178){\rotatebox{45}{\makebox(0,0)[l]{$t^D_3$}}}
\put(233.333,20){\rotatebox{-25}{\makebox(0,0)[l]{$t^D_2+80$}}}
\thicklines
\put(80,25){\line(0,1){10}}
\put(243.333,25){\line(0,1){10}}
\put(150,165){\line(0,1){10}}
\put(173.333,165){\line(0,1){10}}
\end{picture}
\caption{Channel constraint as a polygon packing problem, $t^A$ standing for time of arrival (berth time) and $t^D$ for time of departure}\label{figChannel}
\end{figure}
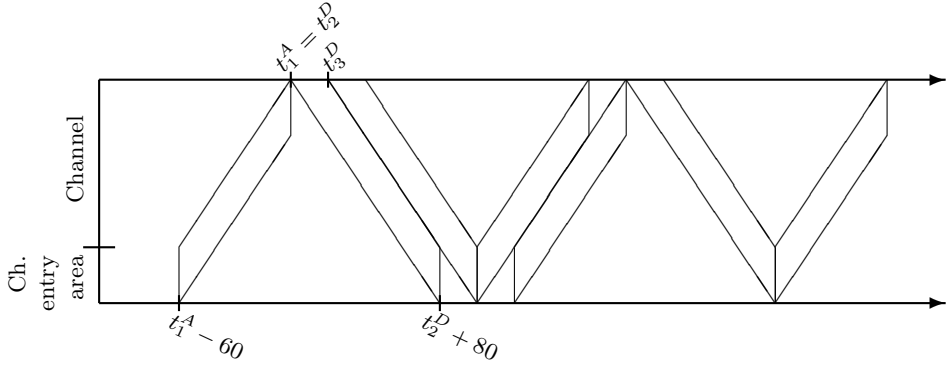

Modeling: the departure time variables $t^\text{Depart}_v$ \eqref{eqTDepart} for all vessels are disconnected from the end of reclaiming by changing the equalities \eqref{eqTDepart} to `greater-or-equal'. Similarly, the arrival time variables $t^\text{Berth}_v$ \eqref{eqBerth} are disconnected from the start of reclaiming, demanding to be earlier or simultaneously. The process of entering or leaving the port for each vessel is partitioned in 20-minute intervals and some of these intervals are demanded to be disjunct among different vessels.

In detail, we represent the departure process by four 20-minute jobs $D_0,D_{20},D_{40},D_{60}$ starting at the following time points:
\be
t^{D_0}_v=t^\text{Depart}_v, \ t^{D_{20}}_v=t^{D_0}_v+20, \ t^{D_{40}}_v=t^{D_0}_v+40, \ t^{D_{60}}_v=t^{D_0}_v+60, \quad\forall v\in V
\ee
Similarly, the arrival process is represented by three 20-minute jobs $A_{60},A_{40},A_{20}$:
\be
t^{A_{60}}_v=t^\text{Berth}_v-60, \ t^{A_{40}}_v=t^\text{Berth}_v-40, \
t^{A_{20}}_v=t^\text{Berth}_v-20,  \quad\forall v\in V
\ee
According to the assumption, the jobs ${D_i}$ of all vessels are all disjoint for a fixed $i\in\{0,20,40,60\}$, as well as jobs $A_j$ for a fixed $j\in\{20,40,60\}$. Moreover, any job $D_i$ is disjoint from any job $A_j$, $\forall i,j$. Let us define vectors $t^{DA}_{ij}\in\Z^{2|V|}$ containing the start times of all jobs $t^{D_i}$ and $t^{A_j}$:
$$
t^{DA}_{ij}=\{t^{D_i}_1,\dots,t^{D_i}_{|V|},t^{A_j}_1,\dots,t^{A_j}_{|V|}\},\qquad i\in\{0,20,40,60\}, \ j\in\{20,40,60\}
$$
  Both restrictions can be modelled by the following 12 \texttt{disjunctive} constraints  \cite{MZN} ensuring  disjointness of the jobs:
\be
\disjunctive(t^{DA}_{ij}, 20\times\mathbf{1}),\qquad i\in\{0,20,40,60\}, \ j\in\{20,40,60\}
\ee
Note that it is possible to reduce the number of variables and constraints, e.g., by disregarding variables $t^{D_{20}}$ and $t^{D_{40}}$ and the constraints involving them, but this showed no significant speed up in the computational experiments.

\paragraph{Flexible stacking volumes:} we can allow some deviation in the daily stacking portions for each stockpile $s\in S$, let's say, the fraction of $\delta^{\tonn}\in[0,1)$ of the nominal daily volume $\tonn^\text{daily}_s$. Thus, we introduce stacking volume variables $v^s_d$ for each of the stacking days $d\in\bigl\{1,\dots,\frac{d^S_s}{1440}\bigr\}$ of each stockpile $s\in S$:
\begin{align}
  v^{\min}_s &= \lfloor(1-\delta^{\tonn})\tonn^\text{daily}_s\rfloor,\qquad&\forall s\in S
  \\
  v^{\max}_s &= \lceil(1+\delta^{\tonn})\tonn^\text{daily}_s\rceil,\qquad&\forall s\in S
  \\
  \stkDays_s &= d^S_s/1440,\qquad &\forall s\in S
  \\
  v^s_1,\ldots,v^s_d&\in\{v^{\min}_s,\dots,v^{\max}_s\}, \qquad d\in\{1,\dots,\stkDays_s\}, \quad &\forall s\in S
  \\
  v^s_1+\dots+v^s_d&=\tonn^\text{daily}_s\stkDays_s,\qquad &\forall s\in S
\end{align}
To formulate a cumulative stacking capacity constraint, we need derived variables for the stacking times:
\be
 t^{SD}_{sd} = t^{S1}_s+d-1, \qquad d\in\{1,\dots,\stkDays_s\}, \quad\forall s\in S
\ee
Deriving stacking volume variables $v^{sk}_d$ for each stream $k\in\{1,2,3\}$ similar to \eqref{eqStStream}, we obtain the cumulative capacity constraints
\begin{gather}
\cumulative(t^{SD}, \mathbf1, v, \tonn^\text{DIT})
\\
\cumulative(t^{SD}, \mathbf1, v^k, \tonn^\text{SS}_k), \qquad k\in\{1,2,3\}
\end{gather}

%

\paragraph{$\mathbf{\pm}$1 day stacking duration:} we make the upper and lower bounds for daily stack volumes variable (${v^{\min}_s}'$, ${v^{\max}_s}'$)\NEW{} and introduce stacking duration variables ${d^S_s}'$ to be used instead of the constants $d^S_s$:
\begin{gather}
  \delta^-_s = \delta^{\tonn}{{\tonn^\text{daily}_s}}
  \frac{\stkDays_s}{\stkDays_s+1},\qquad\forall s\in S
  \\
  \delta^+_s = \delta^{\tonn}{{\tonn^\text{daily}_s}}
  \frac{\stkDays_s}{\stkDays_s-1},\qquad\forall s\in S
  \\
  {v^{\min}_s}' = \Bigl\lfloor{{\tonn^\text{daily}_s}}
  \frac{\stkDays_s}{{d^S_s}'/1440}\Bigr\rfloor - \delta^-_s,\qquad\forall s\in S
  \\
  {v^{\max}_s}' = \Bigl\lceil{{\tonn^\text{daily}_s}}
  \frac{\stkDays_s}{{d^S_s}'/1440}\Bigr\rceil + \delta^+_s,\qquad\forall s\in S
  \\
{d^S_s}' 
    = \begin{cases} {d}^S_i + 1440, &v^s_{\stkDays_s + 1}\ge\max\{1,{v^{\min}_s}'\}, \\
            {d}^S_i, &v^s_{\stkDays_s}\ge\max\{1,{v^{\min}_s}'\} \wedge v^s_{\stkDays_s+1}=0, \\
            {d}^S_i - 1440, &v^s_{\stkDays_s}=v^s_{\stkDays_s+1}=0,
    \end{cases}
    \qquad\forall s\in S
\end{gather}



\paragraph{Dwell post-processing:} Stockpile dwell is the time between
stacking and reclaiming when the stockpile just occupies space on the
pad. Already in \cite{SavSm13} it was found that reducing dwell immediately
after placing a stockpile is generally not a good idea because it reduces
resources for later stockpiles. In our tests, we even found it advantageous
to increase dwell during schedule construction. Thus, we only reduce dwell
by post-processing of complete solutions, changing the stacking start
times and daily volumes of stockpiles for groups of each 5 vessels.

\end{subequations}
}

\subsection{Solver search strategy}\label{secSearchStrat}

\GLEB{As discussed below, it is difficult to obtain good solutions by trying to solve a complete problem instance in a single solver call. Instead, we heuristically decompose the problem into smaller parts by visibility horizons and LNS. However, when solving each smaller part, the search strategy of the CP solver is important.}

Many Constraint Programming models benefit from a custom search
strategy for the solver.  Similar to packing problems
\cite{ClaConstr}, we found it advantageous to separate branching
decisions by groups of variables. We start with the most important
variables --- departure times of all ships (equivalently, delays). \GLEB{This proved helpful to quickly find good solutions}.
Then we fix all reclaim starts, pads, reclaimers, stack starts, and pad
positions. For most of the variables, we use the dichotomous
strategy \texttt{indomain\_split} for value selection, which divides
the current domain of a variable in half and tries first to find a
solution in the lower half. However, pads are assigned randomly, and
reclaimers are assigned preferring lower numbers for odd vessels and
higher numbers for even vessels. Pad positions are preferred so as
to be closer to the native side of the chosen reclaimer, which
corresponds to the idea of opportunity costs in \cite{SavSm13}.
\llong{It appears best not to specify any strategy for daily
stacking volumes, however we make sure that bigger values are tried
first for earlier days. Let us call this
strategy 
\textsc{LayerSearch(1,\ldots,$|V|$)} because we start with all vessels'
departure times, continue with reclaim times, pads, etc. \GLEB{Some experimentation with this strategy is discussed in the results section.}}

In the greedy and LNS heuristics described next, some of the
variables are fixed and the model optimizes only the remaining
variables. For those free variables, we apply the search strategy
described above.

\subsection{\GLEB{Initial feasible solutions: a truncated} search heuristic}\label{secEHRH}

\GLEB{In one CP solver call,} it is difficult to obtain feasible solutions for large
instances in a reasonable amount of time. Moreover, even for average-size
instances, if a feasible solution is found, it is usually bad.
Therefore, we apply a divide-and-conquer strategy which schedules
vessels by groups (e.g., solve vessels 1--5, then vessels 6--10, then
vessels 11--15, etc.). For each group, we allow the solver to
run for a limited amount of time, and,
if feasible solutions are found, take the best of these, or, if no
feasible solution is found, we reduce the number of vessels in the
group and retry. We refer to this scheme as the \emph{extending
horizon} (\EH) heuristic.
%
This heuristic is generalized in Section~\ref{secVH}.



\subsection{Large neighbourhood search}\label{secLNS}

After obtaining a feasible solution, we try to improve it by re-optimizing
subsets of variables while others are fixed to their current values,
a \emph{large neighbourhood search} approach \cite{LNSPis}.
We can apply this improvement approach to both
complete solutions (\emph{global LNS}) or only for the current visibility
horizon (see Section~\ref{secVH}).
The free variables used in the large neighbourhood search
are the decision variables associated with certain stockpiles.
\par
\begin{figure}[tb]
\includegraphics[width=0.47\textwidth]{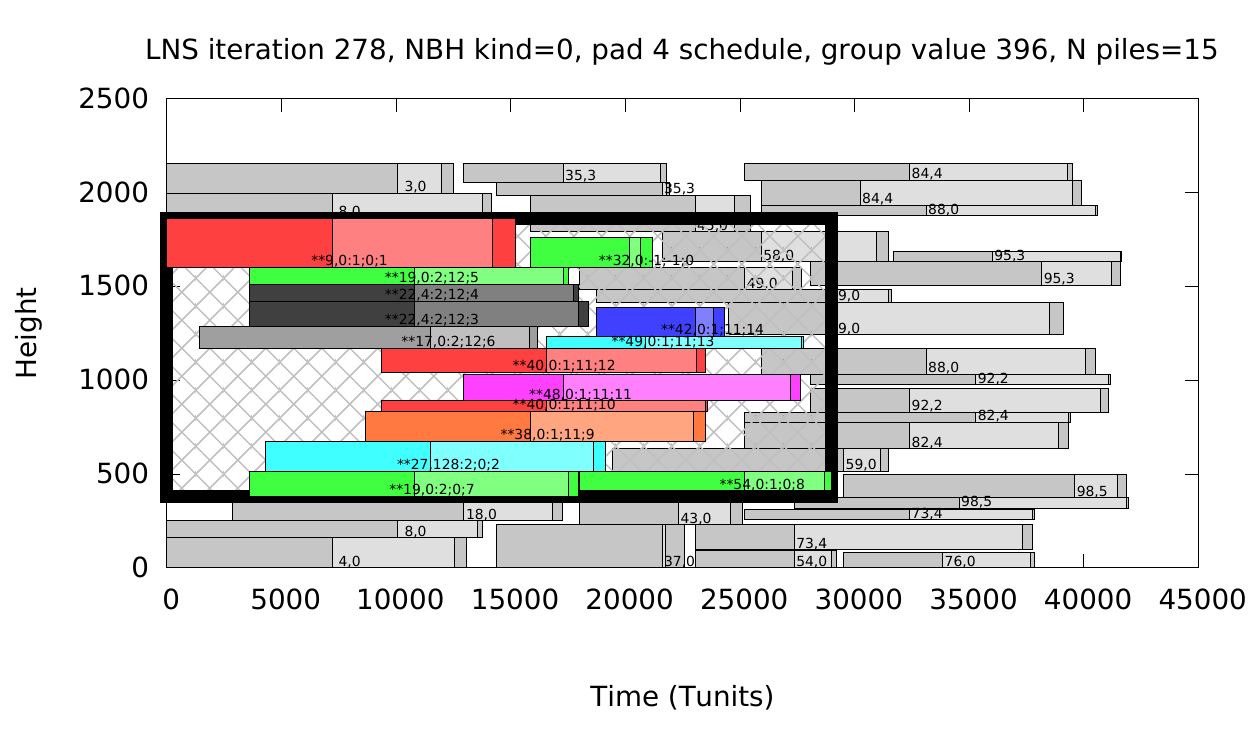}
\hfill
\includegraphics[width=0.47\textwidth]{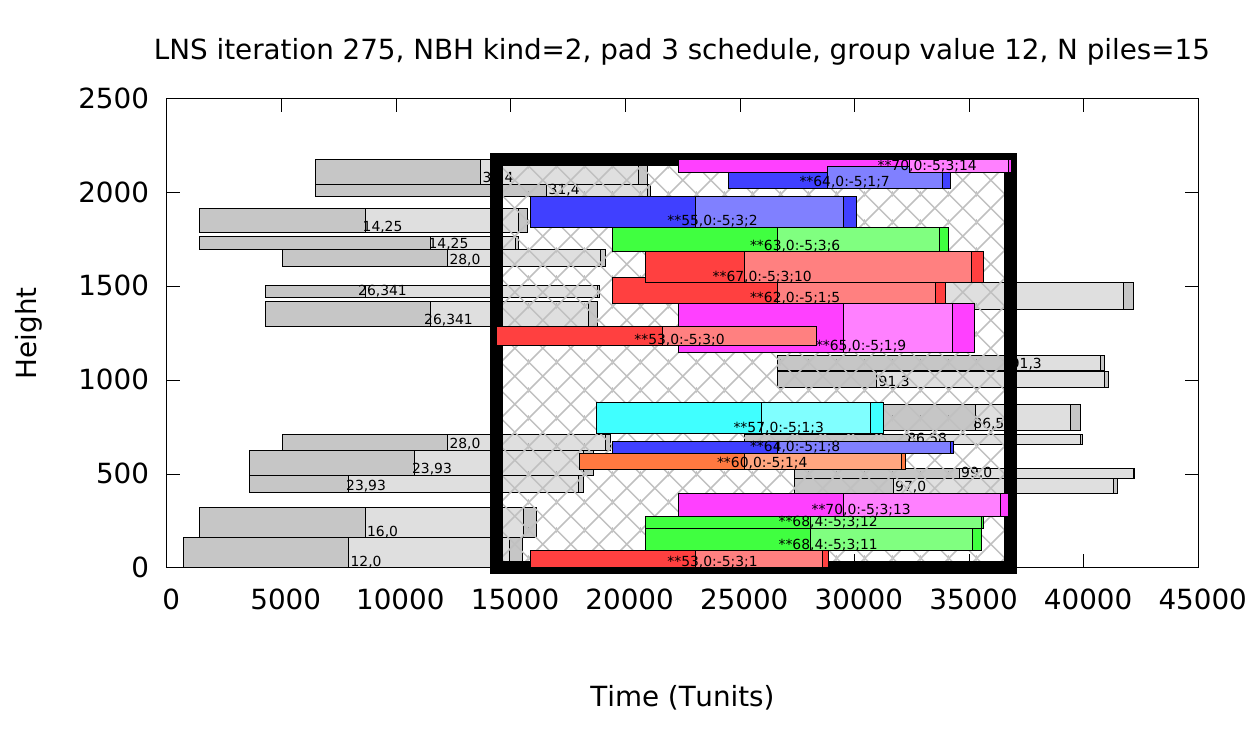}
\caption{Examples of LNS neighbourhoods: spatial (left) and
time-based (right)}\label{figNBH}
\end{figure}

\bigskip\noindent
\textbf{Neighbourhood construction methods.} We
consider a number of methods for choosing which stockpile groups to
re-optimize (the \emph{neighbourhoods}):
\begin{description}\parskip0em\parsep0em\itemsep0em
\item[\quad Spatial] Groups of stockpiles located close to each other on one pad, measured in terms of their space-time location.
\item[\quad Time-based (finish)] Groups of stockpiles on at most two pads with similar reclaim end times.
\item[\quad Time-based (ETA)] Groups of stockpiles on at most two pads belonging to vessels with similar estimated arrival times.
\end{description}
Examples of a spatial and a time-based neighbourhood are given in
Figure~\ref{figNBH}.

First, we randomly decide which of the three types of neighbourhood
to use. Next, we construct all neighbourhoods of the selected type.
Finally, we randomly select one neighbourhood for resolving.

Spatial neighbourhoods are constructed as follows. In order to
obtain many different neighbourhoods, every stockpile seeds a
neighbourhood containing only that stockpile. Then all
neighbourhoods are expanded. Iteratively, for each neighbourhood,
and for each direction right, up, left, and down, independently, we
add the  stockpile on the same pad that is first met by the sweep line going in that direction, after the sweep line has touched the smallest enclosing rectangle of the stockpiles currently in
the neighbourhood. We then add all stockpiles contained in the new
smallest enclosing rectangle. We continue as long as there are
neighbourhoods containing fewer than the target number of stockpiles.

Time-based neighbourhoods are constructed as follows. Stockpiles are
sorted by their reclaim end time or by the ETA of the vessels
they belong to. For each pair of pads, we collect all maximal
stockpile subsequences of the sorted sequence of up to a target length,
with stockpiles allocated to these pads.


Having constructed all neighbourhoods of the chosen type, we
randomly select one neighborhood of the set. The probability of
selecting a given neighborhood is proportional to its
\emph{neighborhood value}: if the \emph{last, but not all}
stockpiles of a vessel is in the neighborhood, then add the vessel's
delay; instead, if \emph{all stockpiles} of a vessel are in the
neighborhood, then add 3 times the vessel's delay.


We denote the iterative large neighbourhood search method by
\textsc{LNS$(k_{\max},n_{\max},\delta)$}, where for at most
$k_{\max}$ iterations, we re-optimize neighborhoods of up to
$n_{\max}$ stockpiles chosen using the principles
outlined above, requiring that the total delay decreases at least by
$\delta$ minutes in each iteration. The objective is again to
minimize the total delay \eqref{eqObj}.


\subsection{Limited visibility horizon}\label{secVH}

In the real world, only a limited number of vessels is known in advance.
We model this as follows: the current \emph{visibility horizon} is
$N$ vessels. We obtain a schedule for the $N$ vessels and fix the
decisions for the first $F$ vessels. Then we schedule vessels
$F+1,\dots,F+N$ (making the next $F$ vessels visible) and so on. Let
us denote this approach by \VH{} $N/F$. Our default visibility
horizon setting is \VH{} 15/5, with the schedule for each visibility
horizon of 15 vessels obtained using \EH{} from
Section~\ref{secEHRH} and then (possibly) improved by LNS(30, 15,
12), i.e., 30 LNS iterations with up to 15 stockpiles in a
neighbourhood, requiring a total delay improvement of at least 12
minutes. We used only time-based neighbourhoods in this case,
because for small horizons, spatial neighbourhoods on one pad are
too small. (Note that the special case \VH{} 5/5 without LNS is equivalent to
the heuristic \EH.)


\section{An adaptive scheme for a heuristic from the literature}\label{secAG}

The truncated tree search (TTS) greedy heuristic \cite{SavSm13} processes vessels
according to a given sequence. It schedules a vessel's stockpiles taking the
vessel's delay into account. It performs a partial lookahead by
considering \emph{opportunity costs} of a stockpile's placement,
which are related to the remaining flexibility of a reclaimer.
However, it does not explicitly take later vessels into account;
thus, the visibility horizon of the heuristic is one vessel. The
heuristic may perform backtracking of its choices if the continuous
reclaim time limit cannot be satisfied.

The default version of TTS processes vessels in their ETA order. We
propose an adaptive framework for this greedy algorithm. This
framework might well be used with the Constraint Programming
heuristic from Section~\ref{secEHRH}, but the latter is slower.
\llong{TTS Greedy does not support the additional constraints from
Section~\ref{secFurtherConstr}, so we compare it with the CP
approach only on the basic model.} Below we present the adaptive
framework, then highlight some modeling differences between CP and
TTS.

\subsection{Two-phase adaptive greedy heuristic (AG)}\label{secAGag}

The TTS greedy heuristic processes vessels in a given order. We
propose an adaptive scheme consisting of two phases. In the first
phase, we iteratively adapt the vessel order, based on vessels'
delays in the generated solutions. In the second phase, earlier
generated orders are randomized. Our motivation to add the
randomization phase was to compare the adaptation principle to pure
randomization.

For the first phase, the idea is to prioritize vessels with large
delays.  We introduce vessels' ``weights'' which are initialized to
the ETAs. In each iteration, the vessels are fed to TTS in order of
non-decreasing weights. Based on the generated solution, the weights
are updated to an average of previous values and ETA minus a
randomized delay; etc. We tried several variants of this principle and
the one that seemed best is shown in Figure~\ref{figAlgAdaptive}, Phase~1.
The variable ``oldWFactor'' is the factor of old weights when averaging them with new values, starting from iteration~1 of Phase~1.

In the second phase, we randomize the orderings obtained in Phase~1.
Each iteration in Phase~1 generated a vessel order
$\mathcal o = (v_1,\dots,v_{|V|})$. Let $\mathcal O = (\mathcal o_1,
\dots, \mathcal o_k)$ be the list of orders generated in Phase~1 in
non-decreasing order of TTS solution value. We select an order with
index $k_0$ from $\mathcal O$ using a truncated geometric
distribution with parameter $p = p_1$, TGD($p$), which has the
following probabilities for indexes $\{1, \dots, k\}$:
\be\nonumber
P[1]=p+(1-p)^k,\quad P[2]=p(p-1), \quad P[3]=p(p-1)^2, \dots, \quad
P[k]=p(p-1)^{k-1} \ee The rationale behind this distribution is to respect the ranking of obtained solutions.
A similar order randomization principle was used, e.g., in
\cite{LeshBLD}.
 Then we modify the selected order $\mathcal o_{k_0}$:
vessels are extracted from it, again using the truncated geometric
distribution with parameter $p = p_2$, and are added to the end of
the new order $\widetilde{\mathcal o}$. Then TTS is executed with
$\widetilde{\mathcal o}$ and $\widetilde{\mathcal o}$ is inserted
into $\mathcal O$ in the position corresponding to its objective
value.
\begin{figure}[tb]\centering
\fbox{
\begin{minipage}{\textwidth}
\begin{tabbing}
\qquad\=\qquad\=\qquad\=\quad\=\kill
Algorithm AG($k_1,k_2$)\\
INPUT: Instance with $V$ set of vessels; $k_1,k_2$ parameters \\
FUNCTION $\rnd(a,b)$ returns a pseudo-random number uniformly distributed in $[a,b)$ \\
Initialize weights: \ $\mathcal W_v=\veta_v$, \ \ $v\in V$ \\
\textbf{for} $k=\overline{0,k_1}$ \qquad\qquad\qquad \qquad\qquad\qquad \qquad\qquad\qquad \qquad\qquad\qquad [PHASE 1]\\
\>Sort vessels by non-decreasing values of $\mathcal W_v$,\\\>\>giving vessels' permutation $\mathcal o=(v_1,\dots,v_{|V|})$\\
\>Run TTS Greedy on $\mathcal o$\\
\>Add $\mathcal o$ to the sorted list $\mathcal O$\\
\>Set $\text{oldWFactor} = \rnd(0.125,1)$ \qquad\qquad\qquad\qquad\qquad // ``Value of history''\\
\>Set $\mathcal D_v$ to be the delay of vessel $v\in V$ \\
\>Let $\mathcal W_v = \text{oldWFactor} \cdot\bigl( \mathcal W_v + (\veta_v - \rnd(0,1) \cdot \mathcal D_v)\bigr), \ \ v \in V$ \\
\textbf{end for}\\
\textbf{for} $k=\overline{1,k_2}$ \qquad\qquad\qquad \qquad\qquad\qquad \qquad\qquad\qquad \qquad\qquad\qquad [PHASE 2]\\
\>Select an ordering $\mathcal o$ from $\mathcal O$ according to TGD($0.5$)\\
\>Create new ordering $\widetilde{\mathcal o}$ from $\mathcal o$,\\
\>\>extracting each new vessel according to TGD($0.85$)\\
\>Run TTS Greedy with the vessel order $\widetilde{\mathcal o}$\\
\>Add the new ordering $\widetilde{\mathcal o}$ to the sorted list $\mathcal O$\\
\textbf{end for}
\end{tabbing}
\end{minipage}
}
\caption{The adaptive scheme for the greedy heuristic.}\label{figAlgAdaptive}
\end{figure}
We denote the
algorithm by AG($k_1,k_2$), where $k_1,k_2$ are the number of
iterations in Phases 1 and 2, respectively. Note that AG($k_1,0$) is
a pure Phase~1 method, while AG($0,k_2$) is a pure randomization
method starting from the ETA order.



\subsection{Differences between the approaches}\label{secDiff}

\short{The model used by both methods is essentially identical, but
there are small technical differences: the CP model uses discrete
time and space and tonnages (minutes, meters, and tons), and
discretizes possible stacking start times to be 12 hours apart. The
discretized stacking start times reduce the search space, and may
diminish solution quality, but seem reasonable given the coarse
granularity of the stacking constraints imposed. The greedy method does
not implement the berth constraints.  If we remove them from the CP
model it is solved slower, but the delay is hardly affected, so we
always include them. }

\llong{
The basic problem options discussed in Sections \ref{secBAsicDescr} and
\ref{secModelBasicCP} are common to both methods. However there are small,
mainly technical differences.

In particular, Constraint Programming works with discrete time and space. In the CP model, we chose the discretization of stacking start times to be 12 hours, which reduces the search space (and thus may diminish solution quality) but is precise enough for the cumulative stacking modeling by streams. All these differences are summarized in Table~\ref{tabDiffMethods}.
 \begin{table}[tb]\centering
 \caption{Differences between the approaches}\label{tabDiffMethods}
 \bigskip
 \begin{tabular}{l|r}
 \makebox[0.3\textwidth]{\textbf{TTS Greedy}} & \makebox[0.3\textwidth]{\textbf{MiniZinc}} \\
   & \\
 Continuous time & Time discretization = 1 min \\
 Continuous position & Position discretization = 1 meter \\
 Continuous stacking volumes \qquad { } & Stacking unit = 1 ton \\
 Continuous stacking start &\qquad Stacking start at 12am and 12pm only 
 \end{tabular}
 \end{table}

Moreover, MiniZinc allows for relatively easy addition of further options to the model, which was discussed in Section~\ref{secFurtherConstr}. As we mentioned in Section~\ref{secBAsicDescr}, in \cite{SavSm13} the authors do not consider the delay of the first and last 8 vessels. In the current paper we optimize over all vessels in both methods.
}

%
\section{Experiments}\label{sec:exp}
%

After describing the experimental set-up, we illustrate the test
data. \GLEB{We start the results presentation with methods to obtain initial solutions. We continue with} the value of
information represented by the visibility horizon. Using the
\llong{basic }model from Section~\ref{secModelBasicCP} we compare
the Constraint Programming approach to the TTS heuristic and the
adaptive scheme from Section~\ref{secAG}. \llong{Then the extended
model options from Section~\ref{secFurtherConstr}
  and other characteristics are tested.
}


The Constraint Programming models in the MiniZinc language were
created by a master program written in C++, which was compiled in
GNU C++.

The adaptive framework for the TTS heuristic and the TTS heuristic
itself were implemented in C++ too. The MiniZinc models were
processed by the finite-domain solver Opturion CPX 1.0.2
\cite{CPXOpturion} which worked single-threaded on an
Intel\textsuperscript{\textregistered}
Core\textsuperscript{\texttrademark} i7-2600 CPU @ 3.40GHz under
Kubuntu~13.04 Linux.

The Lazy Clause Generation~\cite{lazyj} technology seems to be
essential for our approach. Another CP
solver, Gecode 4.3.0 \cite{Gecode}, failed to solve some 1-vessel subproblems, finding no feasible solutions in several hours. Packing problems are highly combinatorial, and this is where learning is the
most advantageous. Moreover, some other LCG solvers than CPX did not work well,
since the solving obviously relies on lazy literal creation.

The solution of a MiniZinc model works in 2 phases. At first, it is
\emph{flattened}, i.e., translated into a simpler language FlatZinc.
Then the actual solver is called on the flattened model. Time limits
were imposed only on the second phase; in particular, we allowed at
most 60 seconds in the \EH{} heuristic and 30 seconds in an LNS
iteration, see
\llong{Section~\ref{secCP}}\short{Section~\ref{secModelBasic}} for
their details. However, reported times contain also the flattening
which took a few seconds per model on average.

\begin{figure}[tb]
\llong{
\includegraphics[page=5,width=\textwidth]{358_delays_3_41_and_7_04}
\caption{Vessel delay profiles in the solutions of the instance 1..358 for
  the basic model and basic + ``stack all before reclaim'' and ``at most 3
  reclaimers active'', obtained with the visibility horizon 15/5.}
}
\short{
\includegraphics[width=\textwidth]{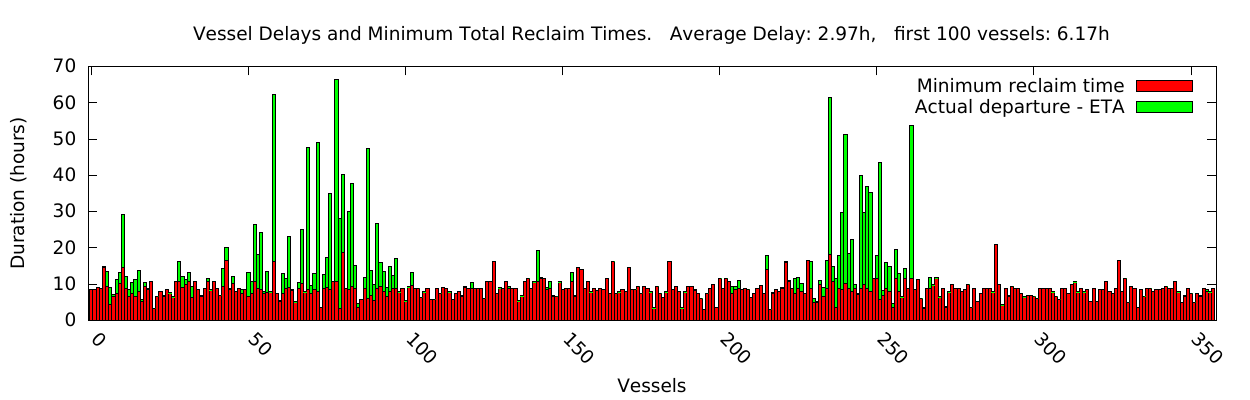}
\caption{Vessel delay profile in a solution of the instance 1..358.}
}
\label{figDelays}
\end{figure}

In \EH{} and LNS, when writing the models with fixed subsets of the
variables, we tried to omit as many irrelevant constraints as
possible. In particular, this helped reduce the flattening time. For
that, we imposed an upper bound of 200 hours on the maximal delay of
any vessel (in the solutions, this bound was never achieved, see
Figure~\ref{figDelays} for example).

The default visibility horizon setting for our experiment, see
Section~\ref{secVH}, is \VH{} 15/5: 15 vessels visible, they are
approximately solved by \EH{} and (possibly) improved by LNS(30, 15,
12); then the first 5 vessels are fixed, etc. Given the above time
limits on an \EH{} or LNS iteration, this takes less than 20 minutes
to process each current visibility horizon and has shown to be
usually much less because many LNS subproblems are proved infeasible
rather quickly. \llong{Only with some extended constraints such as
flexible stacking volumes, \EH{} sometimes needed longer for initial
solutions.}

\llong{
At first we investigate various methods using the basic model from
Section~\ref{secModelBasicCP}.
However for the Constraint Programming
models, it proved computationally advantageous to add the constraint `up to
4 vessels berthed' \eqref{eqMaxBerth}, so we do it always.
}

Our \textbf{test data} is the same as in \cite{SavSm13}. It is historical data with compressed time to put extra pressure on the system. It has the
following key properties:
\begin{itemize}\parskip0em\parsep0em\itemsep0em\topsep0em\topskip0em\partopsep0em
\item 358 vessels in the data file, sorted by their ETAs.
\item One to three stockpiles per vessel, on average 1.4.
\item 
The average interarrival time is 292 minutes.
\item All ETAs are moved so that the first ETA = 10080 (7 days, to
  accommodate the longest build time).
\item Optimizing vessel subsequences of 100 or up to 200 vessels, starting from vessels 1, 21, 41, \dots, 181.
\end{itemize}

Figure~\ref{figDelays} illustrates the test data giving the delay \llong{profiles in solutions}\short{profile in a solution}
for all 358 vessels.
\llong{The two solutions were obtained with the default visibility horizon
  setting \VH{} 15/5. The first solution was obtained with the basic model;
  the second solution also had the additional constraints ``stack all before
  reclaim'' and ``at most 3 reclaimers active''.
}
\short{The solution is obtained with the default visibility horizon
  setting \VH{} 15/5.}
\llong{We see that the average delay is about twice as high in the second solution.}
The most difficult subsequences seem to be the vessel groups 1..100 and
200..270\llong{, whose average delay grows by about the same factor}.
\llong{Below we look at the group 1..100 more closely.}


\subsection{Initial solutions \GLEB{and solver search strategy}}\label{secInit}

First we look at basic methods to obtain schedules for longer sequences of
vessels. This is the \EH{} heuristic from Section~\ref{secEHRH} and TTS
Greedy described in Section~\ref{secAG}, 
which fit into the visibility horizon schemes
\VH{} 5/5 and \VH{} 1/1, respectively. We compare them to an approach to
construct schedules in a single MiniZinc model (method ``ALL'') and to the
standard visibility horizon setting \VH{} 15/5, Section~\ref{secVH}. The
results are given in Table~\ref{tabIni100200} for the 100-vessel and
200-vessel instances.  \GLEB{The results show some interesting properties of the solver search strategy}.

\begin{table}[tb]
\centering
\caption{Solutions for the 100- and up to 200-vessel instances, obtained with \EH, TTS, ALL, and \VH{} 15/5.}
\label{tabIni100200}
\begin{tabular}{r   ||  c       c   |   c       c   |   c       c   |   c       c   ||  c       c   |   c       c   |   c       c   |   c       c   }
    &   \multicolumn{8}{c||}{100 vessels}                                                           &   \multicolumn{8}{c}{Up to 200 vessels}                                                           \\\hline
    &   \multicolumn{2}{c|}{\EH}            &   \multicolumn{2}{c|}{TTS}            &   \multicolumn{2}{c|}{ALL}            &   \multicolumn{2}{c||}{\VH{} 15/5}            &   \multicolumn{2}{c|}{\EH}            &   \multicolumn{2}{c|}{TTS}            &   \multicolumn{2}{c|}{ALL}            &   \multicolumn{2}{c}{\VH{} 15/5}          \\\hline
1st &   obj &   $t$ &   obj &   $t$ &   obj &   $t$ &   obj &   $t$ &   obj &   $t$ &   obj &   $t$ &   obj &   $t$ &   obj &   $t$ \\\hline

1   &   11.77   &   71  &   9.87    &   73  &   13.31   &   275 &   6.17    &   1509    &   6.15    &   170 &   5.09    &   90  &   7.06    &   356 &   3.19    &   1934    \\
21  &   7.01    &   69  &   6.11    &   33  &   9.46    &   275 &   4.19    &   1758    &   3.75    &   142 &   3.25    &   68  &   5.08    &   352 &   2.23    &   2101    \\
41  &   2.54    &   46  &   1.68    &   12  &   2.93    &   271 &   1.31    &   702 &   2.02    &   175 &   1.60    &   62  &   2.62    &   348 &   1.26    &   1465    \\
61  &   0.64    &   42  &   0.61    &   18  &   0.98    &   273 &   0.51    &   214 &   3.59    &   252 &   3.25    &   60  &   5.39    &   351 &   2.63    &   1719    \\
81  &   0.46    &   35  &   0.39    &   18  &   0.54    &   272 &   0.32    &   236 &   3.81    &   139 &   3.40    &   310 &   5.73    &   352 &   2.71    &   2084    \\
101 &   0.33    &   29  &   0.23    &   7   &   0.52    &   272 &   0.19    &   202 &   3.39    &   140 &   3.21    &   62  &   5.14    &   352 &   1.91    &   2444    \\
121 &   0.40    &   27  &   0.38    &   8   &   0.54    &   272 &   0.26    &   169 &   4.79    &   108 &   4.23    &   46  &   4.45    &   360 &   2.33    &   1815    \\
141 &   2.82    &   154 &   1.44    &   20  &   2.59    &   273 &   1.35    &   612 &   4.72    &   220 &   3.26    &   47  &   4.45    &   353 &   2.45    &   2184    \\
161 &   5.13    &   43  &   5.26    &   11  &   7.68    &   273 &   3.67    &   2031    &   3.53    &   101 &   3.26    &   42  &   5.15    &   352 &   2.25    &   2350    \\
181 &   5.45    &   35  &   5.16    &   10  &   8.23    &   273 &   3.84    &   1438    &   3.13    &   70  &   2.93    &   33  &   4.72    &   328 &   2.17    &   1519    \\\hline

Mean    &   3.65    &   55  &   3.11    &   21  &   4.68    &   273 &   2.18    &   887 &   3.89    &   152 &   3.35    &   82  &   4.98    &   350 &   2.31    &   1961
\end{tabular}
\end{table}

Method ``ALL'', obtaining feasible solutions for the whole 100-vessel and
200-vessel instances in a single run of the solver, became possible after a
modification of the default search strategy from
Section~\ref{secSearchStrat}. This did not produce better results however,
so we present its results only as a motivation for iterative methods for
initial construction and improvement.


The default solver search strategy \textsc{LayerSearch(1,\ldots,$|V|$)} proved best for the iterative methods \EH{} and LNS. But feasible solutions of complete instances \emph{in a
  single model} only appeared possible with a modification.  The alternative
strategy can be expressed as
$$
\textsc{GroupLayerSearch(1,\ldots,\text{$|V|$})} = (\textsc{LayerSearch(1,\ldots,5)}; \textsc{LayerSearch(6,\ldots,10)};\ \ldots)
$$
which means: search for departure times of vessels $1,\ldots,5$; then for the reclaim times of their stockpiles; then for their pad numbers; \ldots departure times of vessels $6,\ldots,10$; etc. It is similar to the iterative heuristic \EH{} with the difference that the solver has the complete model and (presumably) takes the first found feasible solution for every 5 vessels.

We had to increase the time limit per solver call: 4 minutes. But the flattening phase took longer than finding a first solution (there are a quadratic number of constraints). Feasible solutions were found in about 1--2 minutes after flattening. We also tried running the solver for longer but this did not lead to better results: the solver enumerates near the leaves of the search tree, which is not efficient in this case. Switching to the solver's default strategy after 300 seconds (search annotation \texttt{cpx\_warm\_start} \cite{CPXOpturion}) gave better solutions, comparable with the \EH{} heuristic.

In Table~\ref{tabIni100200} we see that the solutions obtained by the ``ALL'' method are inferior to \EH.
Thus, for all further tests we used strategy
\textsc{LayerSearch(1,\ldots,$|V|$)} from
Section~\ref{secSearchStrat}. Further, \EH{} is inferior to TTS, both in
quality and running time. This proves the efficiency of the opportunity
costs in TTS and suggests using TTS for initial solutions. However, TTS runs
on original real-valued data and we could not use its solutions in LNS
because the latter works on rounded data which usually has small constraint
violations for TTS solutions. A workaround would be to use the rounded data
in TTS but given \llong{the extended model implemented in MiniZinc only and }the
majority of running time spent in LNS, and for simplicity
we stayed with \EH{} to obtain
starting solutions. The results for \VH{} 15/5 where LNS worked on every
visibility horizon, support this choice.

\subsection{Visibility horizons}\label{secExpVH}

In this subsection, we look at the impact of varying the visibility
horizon settings (Section~\ref{secVH}), including the complete
horizon (all vessels visible). More specifically, we compare $N = 1,
4, 6, 10, 15, 25$, or $\infty$ visible vessels and various numbers
$F$ of vessels to be fixed after the current horizon is scheduled.
For $N = \infty$, we can apply a \emph{global} solution method.
Using Constraint Programming, we obtain an initial solution and try
to improve it by LNS, denoted by \emph{global LNS}, because it
operates on the whole instance. Using Adaptive Greedy
(Section~\ref{secAG}), we also operate on complete schedules.

\begin{figure}[tb]
\includegraphics[page=20,width=0.47\textwidth]{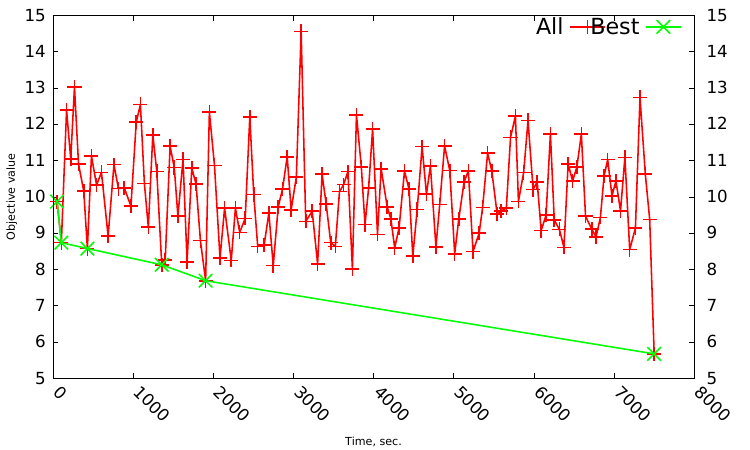}
\hfill
\includegraphics[page=20,width=0.47\textwidth]{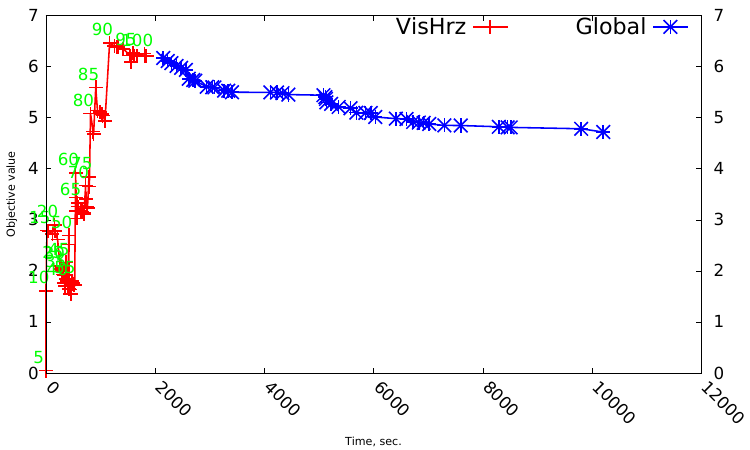}
\caption{Progress of the objective value in AG(130,0) (left) and \VH{} 15/5 + LNS(500, 15, 12) (right), vessels 1..100}\label{figGraphProgress}
\end{figure}

To illustrate the behaviour of global methods, we pick the difficult
instance with vessels 1..100, cf.\ Figure~\ref{figDelays}. A
graphical illustration of the progress over time of the global
methods AG(130,0) and \VH{} 15/5 + LNS(500, 15, 12) is given in
Figure~\ref{figGraphProgress}.

To investigate the value of various visibility horizons, for limited
horizons, we applied the same settings as the standard one
(Section~\ref{secVH}): an initial schedule for the current horizon
is obtained with \EH{} and then improved with LNS(30, 15, 12).
\llong{ Results for the instance 1..100 are given in
Table~\ref{tabVHFirst100}. We see that the best of the limited
visibility horizon settings is \VH{} 25/5; only the global approach,
spending more time, achieves a (much) better solution. }
Table~\ref{tabVH100} gives the average results for all 100-vessel
instances. On average, the global Constraint Programming approach
(500 LNS iterations) gives the best results, but \VH{} 25/5 is
close. Moreover, the setting \VH{} 15/1 which invests significant
effort by fixing only one vessel in a horizon, is slightly better
than \VH{} 25/12, which shows that with a smaller horizon, more
computational effort can be fruitful.

\begin{table}[tb]
\tabcolsep0.4em
\centering
\llong{
 \caption{Visibility horizon trade-off: vessels 1..100, basic model}\label{tabVHFirst100}
\begin{tabular}{r   ||  c       c       c       c       c       c       c       c       c   }
N/F &   1/1 &   4/2 &   6/3 &   10/5    &   15/7    &   15/1    &   25/12   &   25/5    &\bf    15/5+GLNS   \\\hline
Delay, h    &   16.02   &   11.1    &   9.4 &   7.93    &   6.85    &   6.29    &   6.02    &   5.44    &\bf    4.72    \\
\%$\Delta$  &   239\%   &   135\%   &   99\%    &   68\%    &   45\%    &   33\%    &   28\%    &   15\%    &       \\\hline
Time, s &   194 &   430 &   481 &   612 &   2244    &   6542    &   3766    &   7507    &\bf    10204   \\
\%$\Delta$  &   -98\%   &   -96\%   &   -95\%   &   -94\%   &   -78\%   &   -36\%   &   -63\%   &   -26\%   &       \\
\end{tabular}

\par\bigskip                                    }
 \caption{Visibility horizon trade-off: all 100-vessel instances\llong{, basic model}}\label{tabVH100}

\begin{tabular}{r   ||      c   c       c       c       c       c       c       c       c   }
N/F &   1/1 &   4/2 &   6/3 &   10/5    &   15/7    &   15/1    &   25/12   &   25/5    &\bf    15/5+GLNS   \\\hline
Delay, h    &   5.36    &   3.51    &   3.16    &   2.55    &   2.28    &   2.16    &   2.29    &   1.96    &\bf    1.73    \\
\%$\Delta$  &   210\%   &   103\%   &   83\%    &   48\%    &   32\%    &   25\%    &   33\%    &   13\%    &       \\\hline
Time, s &   114 &   188 &   202 &   267 &   916 &   2896    &   1823    &   3354    &\bf    4236    \\
\%$\Delta$  &   -97\%   &   -96\%   &   -95\%   &   -94\%   &   -78\%   &   -32\%   &   -57\%   &   -21\%   &       \\
\end{tabular}
\end{table}

The visibility horizon setting 1/1 produces the worst solutions. The
TTS heuristic of Section~\ref{secAG} also uses this visibility
horizon, but produces better results, see \llong{Tables
\ref{tabCmp100} and \ref{tabCmp200}}\short{Table~\ref{tabCmp100}}.
The reason is probably the more sophisticated search strategy in
TTS, which minimizes `opportunity costs' related to reclaimer
flexibility. At present, it is impossible to implement this complex
search strategy in MiniZinc, the search sublanguage would need
significant extension to do so.

\llong{
Table~\ref{tabVH200} reports similar results for the up to 200-vessel instances, additionally varying the number of fixed vessels in each horizon. Again, \VH{} 25/5 is the best of tested limited horizon configurations, seconded by \VH{} 15/1.
\begin{table}[tb]
 \caption{Visibility horizon trade-off:
averages over all the up to 200-vessel instances, basic model}\label{tabVH200}
\centering      \tabcolsep0.6em
\begin{tabular}{r|      c       c       c       c       c       c       c}
N/F &   1/1 &   4/2 &   6/3 &   10/5    &   15/7    &   25/12   &   15/1    \\\hline
Delay, h    &   5.77    &   3.49    &   3.25    &   2.66    &   2.30    &   2.16    &   2.15    \\
\%$\Delta$  &   194\%   &   78\%    &   66\%    &   36\%    &   17\%    &   10\%    &   10\%    \\\hline
Time, s &   337 &   518 &   518 &   654 &   1739    &   3656    &   6608    \\
\%$\Delta$  &   -97\%   &   -96\%   &   -96\%   &   -95\%   &   -87\%   &   -73\%   &   -51\%   \\
\\
N/F &       &   4/1 &   6/2 &   10/3    &   15/4    &   25/5    &\bf    25/5 + GLNS 300 \\\hline
Delay, h    &       &   3.36    &   3.23    &   2.60    &   2.28    &   2.00    &\bf    1.96    \\
\%$\Delta$  &       &   71\%    &   65\%    &   33\%    &   16\%    &   2\% &       \\\hline
Time, s &       &   917 &   721 &   1211    &   2442    &   7525    &\bf    13479   \\
\%$\Delta$  &       &   -93\%   &   -95\%   &   -91\%   &   -82\%   &   -44\%   &       \\
\end{tabular}
\end{table}
}

\subsection{Comparison of Constraint Programming and Adaptive Greedy}
To compare the Constraint Programming and the AG approaches, we select the
following methods:
\llong{
\defssLMD7.8em
\defssLW8.2em
\begin{defss}
\defssitem{\bf\VH{} 15/5} Visibility horizon 15/5
\defssitem{\bf\VH{} 15/5+G} Visibility horizon 15/5, followed by global LNS 500
\defssitem{\bf\VH{} 25/5} Visibility horizon 25/5
\defssitem{\bf\VH{} 25/5+G} Visibility horizon 25/5, followed by global LNS 300
\defssitem{\bf AG$_1$} TTS Greedy, one iteration on the ETA order
\defssitem{\bf AG$_{500/500}$} Adaptive greedy, 500 iterations in both phases
\defssitem{\bf AG$_{1000/0}$} Adaptive greedy, 1000 iterations in Phase I only
\end{defss}
}
\short{{\bf\VH{} 15/5} Visibility horizon 15/5;
{\bf\VH{} 15/5+G} Visibility horizon 15/5, followed by global LNS 500;
{\bf\VH{} 25/5} Visibility horizon 25/5;
{\bf AG$_1$} TTS Greedy, one iteration on the ETA order;
{\bf AG$_{500/500}$} Adaptive greedy, 500 iterations in both phases;
{\bf AG$_{1000/0}$} Adaptive greedy, 1000 iterations in Phase I only.
}
The results for the 100-vessel instances are in Table~\ref{tabCmp100}\llong{, for
the up to 200-vessel instances in Table~\ref{tabCmp200}}. The pure-random
configuration of the Adaptive Greedy, AG$_{0/1000}$, showed inferior
performance, and its results are not given.
We observe superior performance of LNS on majority of the instances.

As discussed in Section~\ref{secDiff}, CP approach works with discrete measures. We experimented with increasing discretization up to 10 minutes, 10 meters, and 10 tonnes. This slightly improved running times but also impaired objective values by several percent. Still, this might be more robust for real-life solutions.

\begin{table}[tb]
 \caption{\llong{Basic model, }100 vessels: \VH{} and LNS vs.\ (adaptive) greedy}\label{tabCmp100}\centering
 \footnotesize
\begin{tabular}{r   ||  c       c   |   c       c       c   |   c       c   ||  c       c   |   c       c       c   |   c       c       c   }
    &   \multicolumn{7}{c||}{Constraint Programming}                                                    &   \multicolumn{8}{c}{TTS Greedy and Adaptive Greedy}                                                          \\\hline
    &   \multicolumn{2}{c|}{\VH{} 15/5}         &   \multicolumn{3}{c|}{\VH{} 15/5$^*$+G}                   &   \multicolumn{2}{c||}{\VH{} 25/5}            &   \multicolumn{2}{c|}{AG$_1$}         &   \multicolumn{3}{c|}{AG$_{500/500}$}                 &   \multicolumn{3}{c}{AG$_{1000/0}$}                   \\\hline
1st &   obj &   $t$ &   obj &   $t$ &   iter    &   obj &   $t$ &   obj &   $t$ &   obj &   $t$ &   iter    &   obj &   $t$ &   iter    \\\hline

1   &   6.17    &   1509    &   4.72    &   10204   &   338 &   5.44    &   7507    &   9.87    &   73  &   5.67    &   9992    &   130 &   5.67    &   7500    &   130 \\
21  &   4.19    &   1758    &   3.17    &   10529   &   494 &   3.30    &   6626    &   6.11    &   33  &   3.63    &   10433   &   333 &   3.25    &   23051   &   991 \\
41  &   1.31    &   702 &   1.24    &   922 &   0   &   1.24    &   3256    &   1.68    &   12  &   1.00    &   11253   &   667 &   1.02    &   12660   &   903 \\
61  &   0.51    &   214 &   0.50    &   279 &   0   &   0.51    &   912 &   0.61    &   18  &   0.54    &   2713    &   126 &   0.54    &   2769    &   126 \\
81  &   0.32    &   236 &   0.32    &   299 &   0   &   0.32    &   1266    &   0.39    &   18  &   0.34    &   11972   &   696 &   0.36    &   5794    &   344 \\
101 &   0.19    &   202 &   0.19    &   285 &   2   &   0.18    &   943 &   0.23    &   7   &   0.21    &   7764    &   771 &   0.22    &   15  &   1   \\
121 &   0.26    &   169 &   0.26    &   258 &   3   &   0.26    &   971 &   0.38    &   8   &   0.28    &   3589    &   525 &   0.29    &   201 &   28  \\
141 &   1.35    &   612 &   0.73    &   4883    &   469 &   0.90    &   1364    &   1.44    &   20  &   0.80    &   4895    &   255 &   0.76    &   12969   &   652 \\
161 &   3.67    &   2031    &   2.50    &   8172    &   369 &   3.51    &   4241    &   5.26    &   11  &   3.24    &   12574   &   845 &   3.78    &   4166    &   284 \\
181 &   3.84    &   1438    &   3.64    &   6525    &   311 &   3.89    &   6450    &   5.16    &   10  &   3.83    &   6818    &   422 &   3.65    &   13843   &   809 \\\hline

Mean    &   2.18    &   887 &\bf    1.73    &   4236    &   199 &   1.96    &   3354    &   3.11    &   21  &\bf    1.95    &   8200    &   477 &\bf    1.95    &   8297    &   427
\end{tabular}               \par\raggedright\medskip
$^*$ For limited visibility horizons, LNS(20,12,12) was applied
\end{table}

\llong{
\begin{table}[tb]
 \caption{Basic model, up to 200 vessels: \VH{} and LNS vs.\ (adaptive) greedy}\centering
 \footnotesize
 \addtolength{\tabcolsep}{-0.5pt}
\label{tabCmp200}
\begin{tabular}{r   ||  c       c   |   c       c       c   |   c       c   ||  c       c   |   c       c       c   |   c       c       c   }
    &   \multicolumn{7}{c||}{Constraint Programming}                                                    &   \multicolumn{8}{c}{TTS Greedy and Adaptive Greedy}                                                          \\\hline
    &   \multicolumn{2}{c|}{\VH{} 15/5}         &   \multicolumn{3}{c|}{\VH{} 25/5+G}                   &   \multicolumn{2}{c||}{\VH{} 25/5}            &   \multicolumn{2}{c|}{AG$_1$}         &   \multicolumn{3}{c|}{AG$_{500/500}$}                 &   \multicolumn{3}{c}{AG$_{1000/0}$}                   \\\hline
1st &   obj &   $t$ &   obj &   $t$ &   iter    &   obj &   $t$ &   obj &   $t$ &   obj &   $t$ &   iter    &   obj &   $t$ &   iter    \\\hline

1   &   3.19    &   1934    &   2.69    &   18246   &   274 &   2.82    &   10442   &   5.09    &   90  &   3.63    &   72884   &   539 &   3.68    &   102290  &   793 \\
    &   2.23    &   2101    &   2.05    &   15671   &   243 &   1.79    &   8573    &   3.25    &   68  &   1.96    &   6496    &   81  &   1.92    &   72050   &   866 \\
    &   1.26    &   1465    &   1.19    &   11681   &   286 &   1.15    &   6369    &   1.60    &   62  &   1.07    &   7461    &   86  &   1.06    &   7580    &   86  \\
.   &   2.63    &   1719    &   2.37    &   14884   &   245 &   1.82    &   6728    &   3.25    &   60  &   2.09    &   38128   &   532 &   1.94    &   32747   &   442 \\
.   &   2.71    &   2084    &   1.87    &   16147   &   291 &   2.14    &   8113    &   3.40    &   310 &   2.21    &   31628   &   252 &   2.04    &   86931   &   663 \\
.   &   1.91    &   2444    &   1.78    &   12404   &   173 &   2.23    &   7842    &   3.21    &   62  &   2.04    &   102010  &   860 &   2.22    &   56555   &   590 \\
    &   2.33    &   1815    &   1.84    &   14295   &   272 &   1.71    &   6380    &   4.23    &   46  &   2.08    &   82210   &   894 &   2.32    &   26512   &   272 \\
    &   2.45    &   2184    &   1.90    &   12304   &   214 &   2.10    &   7361    &   3.26    &   47  &   2.17    &   57381   &   674 &   2.15    &   5451    &   61  \\
    &   2.25    &   2350    &   1.87    &   12729   &   279 &   2.02    &   6638    &   3.26    &   42  &   2.12    &   33929   &   522 &   2.46    &   18109   &   287 \\
181 &   2.17    &   1519    &   2.04    &   6429    &   12  &   2.20    &   6801    &   2.93    &   33  &   1.99    &   16687   &   488 &   1.99    &   16668   &   488 \\\hline

Mean    &   2.31    &   1961    &\bf    1.96    &   13479   &   229 &   2.00    &   7525    &   3.35    &   82  &\bf    2.14    &   44881   &   493 &   2.18    &   42489   &   455
\end{tabular}           
\end{table}
}


\llong{
\subsection{Extended model}\label{secExpExtModel}

This subsection tests the extended model options from Section~\ref{secFurtherConstr} with the default visibility horizon setting, Section~\ref{secVH}.

A sensitivity analysis for the extended constraints is given in
Table~\ref{tabSensNew}:
\begin{itemize}
\item
Allowing different pads for the stockpiles of a same vessel deteriorates the
solutions, similar to the results in \cite{SavSm13}. Our experiments with
the search strategy did not change this conclusion, thus we accept the
``same pad'' strategy as default.
\item 
Stacking all stockpiles for a vessel before reclaiming slightly deterioated
solutions, and adds some cost in producing solutions.
\item 
The constraint of at most 3 reclaimers working simultaneously is very
strong, doubling the delays, and substantially increasing solving time.
\item 
Setting tidal weight of 80,000t (about 20\% of the vessels) has almost no
effect. But even 70,000t, exceeded by about 50\% of the vessels, increases
average delay only a little. The reason is that we update the earliest
departure for delay calculation to the next tidal window border. Without
doing this, the delays are much larger.  The reduced choice caused by tidal
constraints substantially improves solving time.
\item 
The channel constraints are complex to model and thus costly in solving time,
and do significantly extend average delay.
\item 
Allowing 80\% daily stacking volume deviation gives only a small
improvement.
Allowing $\pm1$ day stacking duration reduces the delays by a factor of 3 in
the basic model; but similar schedules arise if we simply
shorten all stacking durations by 1 day, which is computationally simpler.
\item The combined effect of the new constraints ``stack all before reclaim'',
``tidal weight 70,000t'', ``at most 3 reclaimers'', and ``channel constraints''
is nearly the sum of their effects.  
\item 
In this ``full'' setting, reducing the stacking duration by one day has a smaller effect, but still leads to a reduction in delay of 11\%.
\end{itemize}

Table~\ref{tabVHNew} gives visibility horizon trage-offs with the new constraints. The percentage impacts are slightly smaller than in the basic model, Table~\ref{tabVH200}.

As discussed in the end of Section~\ref{secBAsicDescr}, some of the first and last vessels of an instance
are possibly easy to schedule because they have more resources. To investigate the
hardness of processing `middle' vessels, we report average delays excluding 40
`warm-up' and 20 `cool-down' vessels, however in solutions where still the total delay was minimized.
The corresponding visibility horizon trade-offs with the new constraints are given in Table~\ref{tabVHNewMiddle}. As expected, the delay values are larger, cf.\ Table~\ref{tabVHNew}. 
 The effects of the visibility horizon changes are
similar to the original objective.



\begin{table}[tb]
 \caption{Sensitivity for the new constraints}\label{tabSensNew}\centering
\tabcolsep0.5em
\begin{tabular}{rr  ||  c       c   |   c       c   ||  c       c   |   c       c   }

    &       &   \multicolumn{4}{c||}{100 vessels}                           &   \multicolumn{4}{c}{Up to 200 vessels}                           \\\hline
    &       &   obj &   $t$ &   \%$\Delta$obj   &   \%$\Delta t$    &   obj &   $t$ &   \%$\Delta$obj   &   \%$\Delta t$    \\\hline

0)  &   Basic model$^*$ &   2.18    &   887 &       &       &   2.31    &   1961    &       &       \\
    &   Diff pads   &   2.52    &   672 &   16\%    &   -24\%   &   2.52    &   1398    &   9\% &   -29\%   \\
1)  &   Stack all before    &   2.43    &   1018    &   11\%    &   15\%    &   2.43    &   2091    &   5\% &   7\% \\
2)  &   3 reclaimers    &   4.78    &   1279    &   119\%   &   44\%    &   4.41    &   2612    &   91\%    &   33\%    \\
    &   Tidal 80,000t   &   2.26    &   666 &   4\% &   -25\%   &   2.21    &   1657    &   -4\%    &   -16\%   \\
3)  &   Tidal 70,000t   &   2.38    &   705 &   9\% &   -21\%   &   2.38    &   1693    &   3\% &   -14\%   \\
4)  &   Channel &   3.14    &   1519    &   44\%    &   71\%    &   3.07    &   3303    &   33\%    &   68\%    \\
\multicolumn{2}{r||}{80\% daily stk flex}           &   2.22    &   1171    &   2\% &   32\%    &   2.27    &   3108    &   -2\%    &   58\%    \\
    &   $\pm$1 day stk dur  &   0.82    &   884 &   -62\%   &   0\% &   0.74    &   2679    &   -68\%   &   37\%    \\
5)  &   $-$1 day stk dur    &   0.81    &   444 &   -63\%   &   -50\%   &   0.80    &   1164    &   -65\%   &   -41\%   \\\hline
    &   0-4)    &   5.80    &   1472    &   166\%   &   66\%    &   5.38    &   3343    &   133\%   &   70\%    \\
    &   0-5)    &   5.14    &   1418    &   -11\%   &   -4\%    &   4.81    &   3120    &   -11\%   &   -7\%

\end{tabular}
            \par\medskip        \raggedright
$^*$ The constraint on 4 berthed ships is already included by default
\end{table}

\begin{table}[tb]
 \caption{Up to 200: VH trade-off with the new constraints: stack all stockpiles before reclaim; max 3 reclaimers; tidal ship weight 70'000t; channel constraint}\label{tabVHNew}\centering
 \tabcolsep0.5em
\begin{tabular}{r|      c       c       c       c       c       c       c   }
N/F &   1/1 &   4/2 &   6/3 &   10/5    &   15/7    &   25/12   &       \\\hline
Delay, h    &   11.62   &   7.97    &   6.72    &   5.82    &   5.35    &   5.29    &       \\
\%$\Delta$  &   141\%   &   65\%    &   39\%    &   21\%    &   11\%    &   10\%    &       \\\hline
Time, s &   419 &   852 &   933 &   1425    &   2708    &   5123    &       \\
\%$\Delta$  &   -97\%   &   -94\%   &   -94\%   &   -90\%   &   -81\%   &   -65\%   &       \\
\\
N/F &       &   4/1 &   6/2 &   10/3    &   15/4    &   25/5    &\bf    '25/5 + GLNS 300    \\\hline
Delay, h    &       &   7.42    &   6.36    &   5.77    &   5.07    &   4.97    &\bf    4.82    \\
\%$\Delta$  &       &   54\%    &   32\%    &   20\%    &   5\% &   3\% &       \\\hline
Time, s &       &   1517    &   1271    &   1910    &   3580    &   10280   &\bf    14511   \\
\%$\Delta$  &       &   -90\%   &   -91\%   &   -87\%   &   -75\%   &   -29\%   &       \\
\end{tabular}

\par\bigskip
 \caption{Up to 200: VH trade-off with the new constraints, average delay reported under exclusion of the first 40 and last 20 vessels}\label{tabVHNewMiddle}\centering
\begin{tabular}{r|      c       c       c       c       c       c       c   }
N/F &   1/1 &   4/1 &   6/2 &   10/3    &   15/4    &   25/5    &\bf    '25/5 + GLNS 300    \\\hline
Delay, h    &   13.39   &   8.59    &   7.30    &   6.66    &   5.69    &   5.48    &\bf    5.36    \\
\%$\Delta$  &   150\%   &   60\%    &   36\%    &   24\%    &   6\% &   2\% &
\end{tabular}
\end{table}

 \subsection{Dwell reduction.}
As already noticed in \cite{SavSm13}, it is better not to reduce dwell during schedule construction. We found it better to \emph{increase dwell}, in order to free resources for future vessels.  Moreover, our experiments with dwell reduction for limited visibility horizons led to worse results as well. 

Thus, at the moment we only reduce dwell by post-processing the complete
solution. Reclaim times are fixed and we simply try to move the stacking
start time of stockpiles forward as much as possible without violating
stacking capacity constraints. 
This is done for the stockpiles of consecutive groups of 5 vessels.
Average results before and after dwell post-processing are given in
Table~\ref{tabDwell}. Clearly the post-processing massively reduces dwell.
\begin{table}[tb]
\caption{Average dwell before and after dwell reduction, hours}\label{tabDwell}
\centering\tabcolsep1.5em
\begin{tabular}{r|      c       c   |   c       c   |   c       c   }
    &   \multicolumn{4}{c|}{Basic model}                            &   \multicolumn{2}{c}{Basic + 1,2,3,4)}            \\\hline
    &   \multicolumn{2}{c|}{100 vessels}            &   \multicolumn{4}{c}{Up to 200 vessels}                           \\\hline
Mean dwell  &   81  &   6.3 &   78  &   6.6 &   80  &   7.9
\end{tabular}
\end{table}
An example of the effect of dwell reduction on a pad schedule is given in Figure~\ref{figDwell}.
\begin{figure}[tb]
\includegraphics[page=8,width=0.47\textwidth
]{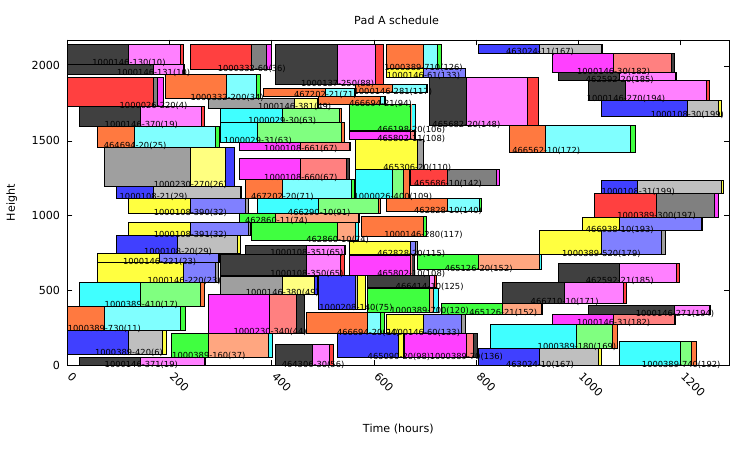}
\hfill
\includegraphics[page=8,width=0.47\textwidth
]{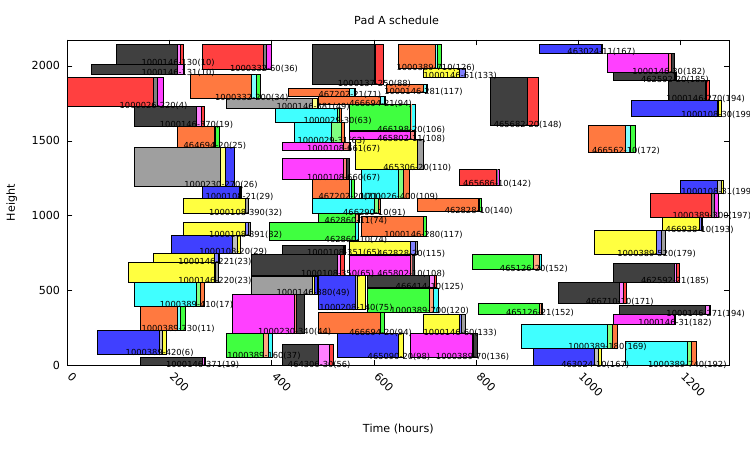}
\caption{Example pad schedules: before and after dwell reduction}\label{figDwell}
\end{figure}




\subsection{Queue jumps}

  In practice, understanding the changes to the vessel processing order, as compared to the ETA order, is important: a customer which has timely arrived to the port entrance is not happy to see too many other vessels arriving later and being processed before him. We did not implement any constraints on such `queue jumps' but observed them in the obtained solutions. Namely, we registered the following values:
\defssLMD5.8em
\defssLW6.2em
\begin{itemize}
\item
 maximal jump (number of vessels) of a vessel behind its ETA-based position, in the vessel ordering by berth start time (arrival)
\item
 same as above, in the ordering by berth end time (departure)
\end{itemize}
Table~\ref{tabFurtherParams} presents the above characteristics of the CP
(\VH{}~15/5, basic model) and AG$_{500/500}$ solutions. We see that jumps by
departure times are a little higher. CP solutions have larger order
permutations.\NEW{} They increase even more after global LNS.

We may consider adding variables and constraints to directly model and
constrain queue jumping in the future. 

\begin{table}[tb]
\caption{Order jumps by arrival, departure time}
\label{tabFurtherParams}\centering
\small
\tabcolsep1.2ex
\begin{tabular}{r   ||  c       c   |   c       c   ||  c       c   |   c       c     }
    &   \multicolumn{4}{c||}{100 vessels}                           &   \multicolumn{4}{c}{Up to 200 vessels}                                                         \\\hline
1st &   \multicolumn{2}{c|}{CP}         &   \multicolumn{2}{c||}{AG}            &   \multicolumn{2}{c|}{CP}         &   \multicolumn{2}{c}{AG}            \\
                                                                                                    \hline
1   &   8   &   11  &   13  &   13  &   14  &   15  &   7   &   9    \\
21    &   11  &   11  &   7   &   10  &   13  &   14  &   7   &   7     \\
41    &   5   &   5   &   1   &   4   &   4   &   8   &   2   &   4   \\
61   &   2   &   4   &   3   &   5   &   11  &   12  &   8   &   9   \\
81   &   1   &   4   &   1   &   4   &   12  &   16  &   7   &   8    \\
101  &   1   &   1   &   1   &   2   &   14  &   17  &   5   &   7   \\
121   &   2   &   4   &   2   &   4   &   16  &   17  &   10  &   10   \\
141   &   5   &   7   &   2   &   4   &   16  &   17  &   5   &   5     \\
161   &   12  &   13  &   6   &   6   &   13  &   14  &   8   &   8   \\
181 &   10  &   10  &   8   &   10  &   11  &   11  &   5   &   5   \\
                                                                                                    \hline
mean    &   5.7 &   7.0 &   4.4 &   6.2 &   12.4    &   14.1    &   6.4 &   7.2 

\end{tabular}           
\end{table}
}

%
\section{Conclusions\llong{ and outlook}}\label{sec:conc}
%

We consider a complex problem involving scheduling and allocation of
cargo assembly in a stockyard, loading of cargoes onto vessels, and
vessel scheduling. We designed a Constraint Programming (CP)
approach to construct feasible solutions and improve them by Large
Neighbourhood Search (LNS).

Investigation of various visibility horizon settings has shown that larger
numbers of known arriving vessels lead to better results. In particular, the
visibility horizon of 25 vessels provides solutions close to the
best found.
\llong{Among the
\GLEB{extended model options}, the strongest impacts arise from: 
the restriction to 3
reclaimers working at any time, which increased the average delay by 91\%;
and the possibility to speed up stacking by 1 day, which reduced delays by
65\% in the basic model and by 11\% with all the new constraints added.

We observed that the constraints involving reclaimer moving speed,
\eqref{eqReclMov} and \eqref{eqReclClash}, were very difficult for the
solver\GLEB{, given that we model reclaim scheduling at a minute granularity}. Simplifying the model to ignore reclaimer movement time
caused the models to be much easier to solve and allowed 
a larger number of vessels to be considered together in one solver call.  
However, in the obtained solutions more than half of the necessary
travel distance between reclaim jobs was to be covered in zero time,
which is unacceptable.

For the basic variant of the model, t}\short{T}he new approach was
compared to an existing greedy heuristic. The latter works with a
visibility horizon of one vessel and, under this setting, produces
better feasible solutions in less time. The reason is probably the
sophisticated search strategy which cannot be implemented in the
chosen CP approach at the moment. To make the comparison fairer, an
adaptive iterative scheme was proposed for this greedy heuristic,
which resulted in a similar performance to LNS.

Overall the CP approach using visibility horizons and LNS generated
the best overall solutions in less time than the adaptive greedy
approach. A significant advantage of the CP approach is that it is
easy to include additional constraints, \short{which we have done in
work not reported here for space reasons.}\llong{as we have done
here.}

\llong{
For further work, it might be important to incorporate other existing types
of stockyard operation. An estimation of the solutions' closeness to optimum
is another goal. The Constraint Programming approach could benefit from
integration of LNS into solver's search strategy.
}


\subsubsection*{Acknowledgments}

\anonymize{The research presented here is supported by ARC linkage grant LP110200524.}
We would like to thank the strategic planning team at HVCCC for many insightful and helpful suggestions\anonymRepl{, Andreas Schutt for hints on efficient modeling in
MiniZinc, as well as to }{ and }Opturion for providing their version of the CPX
solver under an academic license.

\bibliography{../../bib/jourabb,../../bib/jourabbs,../../bib/cp2011,../../bib/cp2004,../../bib/cp2003,../../bib/cp196x,../../bib/cp197x,../../bib/cp198x,../../bib/cp199x}

\begin{thebibliography}{10}
\providecommand{\url}[1]{\texttt{#1}}
\providecommand{\urlprefix}{URL }

\bibitem{Bay10}
Bay, M., Crama, Y., Langer, Y., Rigo, P.: Space and time allocation in a
  shipyard assembly hall. Annals of Operations Research  179(1),  57--76 (2010)

\bibitem{BeldiceanuGlobal94}
Beldiceanu, N., Contejean, E.: Introducing global constraints in {CHIP}.
  Mathematical and Computer Modelling  20(12),  97--123 (1994)

\bibitem{BBSS14small}
Belov, G., Boland, N., Savelsbergh, M.W.P., Stuckey, P.J.: Local search for a
  cargo assembly planning problem. In: Simonis, H. (ed.) Integration of AI and
  OR Techniques in Constraint Programming, Lecture Notes in Computer Science,
  vol. 8451, pp. 159--175. Springer International Publishing (2014)

\bibitem{Gulcz12}
Boland, N., Gulczynski, D., Savelsbergh, M.: A stockyard planning problem. EURO
  Journal on Transportation and Logistics  1(3),  197--236 (2012)

\bibitem{BSHVCC05}
Boland, N.L., Savelsbergh, M.W.P.: Optimizing the {Hunter Valley Coal Chain}.
  In: Gurnani, H., Mehrotra, A., Ray, S. (eds.) Supply Chain Disruptions, pp.
  275--302. Springer London (2012)

\bibitem{ClaConstr}
Clautiaux, F., Jouglet, A., Carlier, J., Moukrim, A.: A new constraint
  programming approach for the orthogonal packing problem. Computers {\&}
  Operations Research  35(3),  944--959 (2008)

\bibitem{HVCCC}
HVCCC: Hunter valley coal chain --- overview presentation (2014),
  http://www.hvccc.com.au/

\bibitem{LeshBLD}
Lesh, N., Mitzenmacher, M.: {BubbleSearch}: A simple heuristic for improving
  priority-based greedy algorithms. Information Processing Letters  97(4),
  161--169 (2006)

\bibitem{LMV02}
Lodi, A., Martello, S., Vigo, D.: Recent advances on two-dimensional bin
  packing problems. Discrete Applied Mathematics  123(1--3),  379--396 (2002)

\bibitem{MZN}
Marriott, K., Stuckey, P.J.: A {MiniZinc} tutorial (2012),
  http://www.minizinc.org/

\bibitem{LCG}
Ohrimenko, O., Stuckey, P.J., Codish, M.: Propagation = lazy clause generation.
  In: Bessi{\`e}re, C. (ed.) Principles and Practice of Constraint Programming
  --- CP 2007, Lecture Notes in Computer Science, vol. 4741, pp. 544--558.
  Springer Berlin Heidelberg (2007)

\bibitem{lazyj}
Ohrimenko, O., Stuckey, P., Codish, M.: Propagation via lazy clause generation.
  Constraints  14(3),  357--391 (2009)

\bibitem{CPXOpturion}
{Opturion Pty Ltd}: Opturion {CPX} user's guide: version 1.0.2 (2013),
  www.opturion.com

\bibitem{LNSPis}
Pisinger, D., Ropke, S.: Large neighborhood search. In: Gendreau, M., Potvin,
  J.Y. (eds.) Handbook of Metaheuristics, International Series in Operations
  Research \& Management Science, vol. 146, pp. 399--419. Springer US (2010)

\bibitem{SavSm13}
Savelsbergh, M., Smith, O.: Cargo assembly planning. EURO Journal on
  Transportation and Logistics pp. 1--34 (2014)

\bibitem{Gecode}
Schulte, C., Tack, G., Lagerkvist, M.Z.: Modeling and programming with {Gecode}
  (2014), www.gecode.org

\bibitem{CumulExpl11}
Schutt, A., Feydy, T., Stuckey, P.J., Wallace, M.G.: Explaining the cumulative
  propagator. Constraints  16(3),  250--282 (2011)

\bibitem{RCSPmax}
Schutt, A., Feydy, T., Stuckey, P.J., Wallace, M.G.: Solving {RCPSP}/max by
  lazy clause generation. Journal of Scheduling  16(3),  273--289 (2013)

\bibitem{Carpet11}
Schutt, A., Stuckey, P.J., Verden, A.R.: Optimal carpet cutting. In: Lee, J.
  (ed.) Principles and Practice of Constraint Programming --- CP 2011, Lecture
  Notes in Computer Science, vol. 6876, pp. 69--84. Springer Berlin Heidelberg
  (2011)

\bibitem{Singh12}
Singh, G., Sier, D., Ernst, A.T., Gavriliouk, O., Oyston, R., Giles, T.,
  Welgama, P.: A mixed integer programming model for long term capacity
  expansion planning: A case study from the hunter valley coal chain. European
  Journal of Operational Research  220(1),  210--224 (2012)

\bibitem{Thomas13}
Thomas, A., Singh, G., Krishnamoorthy, M., Venkateswaran, J.: Distributed
  optimisation method for multi-resource constrained scheduling in coal supply
  chains. International Journal of Production Research  51(9),  2740--2759
  (2013)

\end{thebibliography}
\bibliographystyle{splncs03}

\end{document}